\documentclass{elsarticle}

\usepackage{geometry,t1enc,amsfonts,latexsym,amsmath,amsthm,verbatim,amssymb,graphicx,mathdots,pstool,color,algorithm,algpseudocode,subfig,wrapfig}

\def\N{{\mathbb N}}

\def\C{{\mathbb C}}

\def\R{{\mathbb R}}
\def\Z{{\mathbb Z}}

\def\dia{{\mathrm{diag}}}

\newcommand{\bb}[1]{\boldsymbol{#1}}

\numberwithin{equation}{section}

\graphicspath{/code}

\begin{document}

\title{ESPRIT for multidimensional general grids.}
\author[add1]{Fredrik Andersson}
\ead{fredrik.andersson@erdw.ethz.ch}
\author[add2]{Marcus Carlsson}
\ead{marcus.carlsson@math.lu.se}

\address[add1]{Institute of Geophysics, ETH-Zurich, Sonneggstrasse 5, 8092 Zurich, Switzerland\\}
\address[add2]{Centre for Mathematical Sciences, Lund University, Box 118, SE-22100, Lund,  Sweden\\}

\begin{abstract}
We present a method for complex frequency estimation in several variables, extending the classical one-dimensional ESPRIT-algorithm, and consider how to work with data sampled on non-standard domains, i.e going beyond multi-rectangles.
\end{abstract}

\maketitle

\section{Introduction}
Given equidistant measurements of a signal $f$ and a model order $K$ we seek to approximate $f$ by an exponential sum
\begin{equation}
  \label{expsep1d}f(x)\approx\sum_{k=1}^K c_k e^{\zeta_k x},\quad c_k,\zeta_k\in\C.
\end{equation}
The complex parameters $\zeta_k$ play the role of frequencies in the case when they lie on the imaginary axis, and the problem is therefore often referred to as \emph{frequency estimation}. Once the parameters $\zeta_k$ are found, the subsequent retrieval of the coefficients $c_k$ is a linear problem which can be solved by the least squares method. In this paper, we are interested in multivariate version of this problem.

The frequency estimation problem goes far back in time. Already in 1795, it was discussed by Prony \cite{Prony}. In the noise free case, i.e. when $f$ already is of the form given by the right hand side of \eqref{expsep1d}, it is possible to exactly recover both the complex frequencies and the corresponding coefficients. The approximation problem is of great importance, and several different approaches and methods have been devoted to solving it. Classical methods include ESPRIT \cite{roy1989esprit}, MUSIC \cite{schmidt1986multiple}, Pisarenkos method \cite{pisarenko1973retrieval}, and the matrix pencil method, see e.g. \cite{hua1990matrix} which also provides a good account on the history of the method and precursors, as well as the recent contribution \cite{potts2013parameter} which finds a unifying framework for these methods. From the perspective of functional analysis, a lot of theoretic results were obtained by Adamjan, Arov and Krein (AAK) \cite{adamjan1971analytic}. These connections were further exploited in \cite{beylkin2005approximation} for the construction of optimal quadrature nodes. The connection between approximation theory and the results by Adamjan, Arov and Krein is discussed in more detail in \cite{andersson2011sparse}.

\begin{figure}
	\centering
	\subfloat[][Real part of $f$.]{\includegraphics[width=0.47\linewidth]{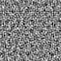}}
	\subfloat[][The corresponding block-Hankel matrix]{\includegraphics[width=0.47\linewidth]{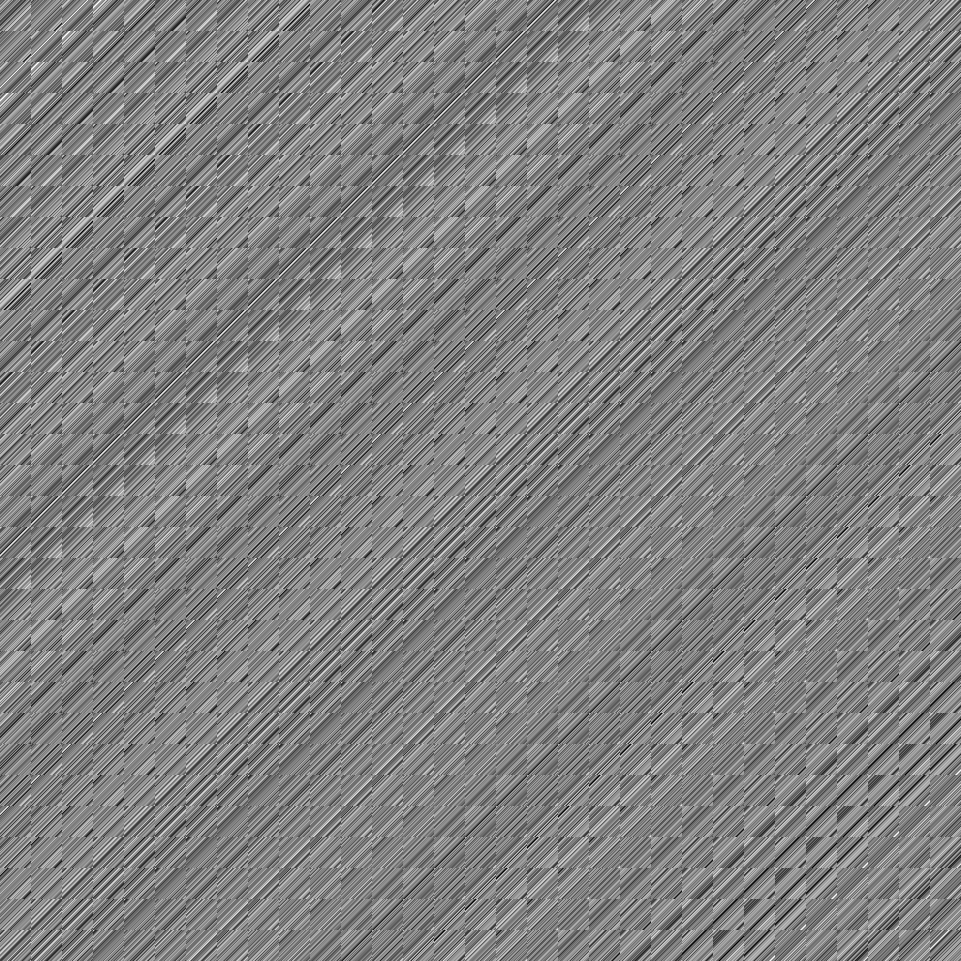}}\\
	\subfloat[][True 2d frequencies in blue dots and estimated frequencies in red circles]{\includegraphics[width=0.47\linewidth]{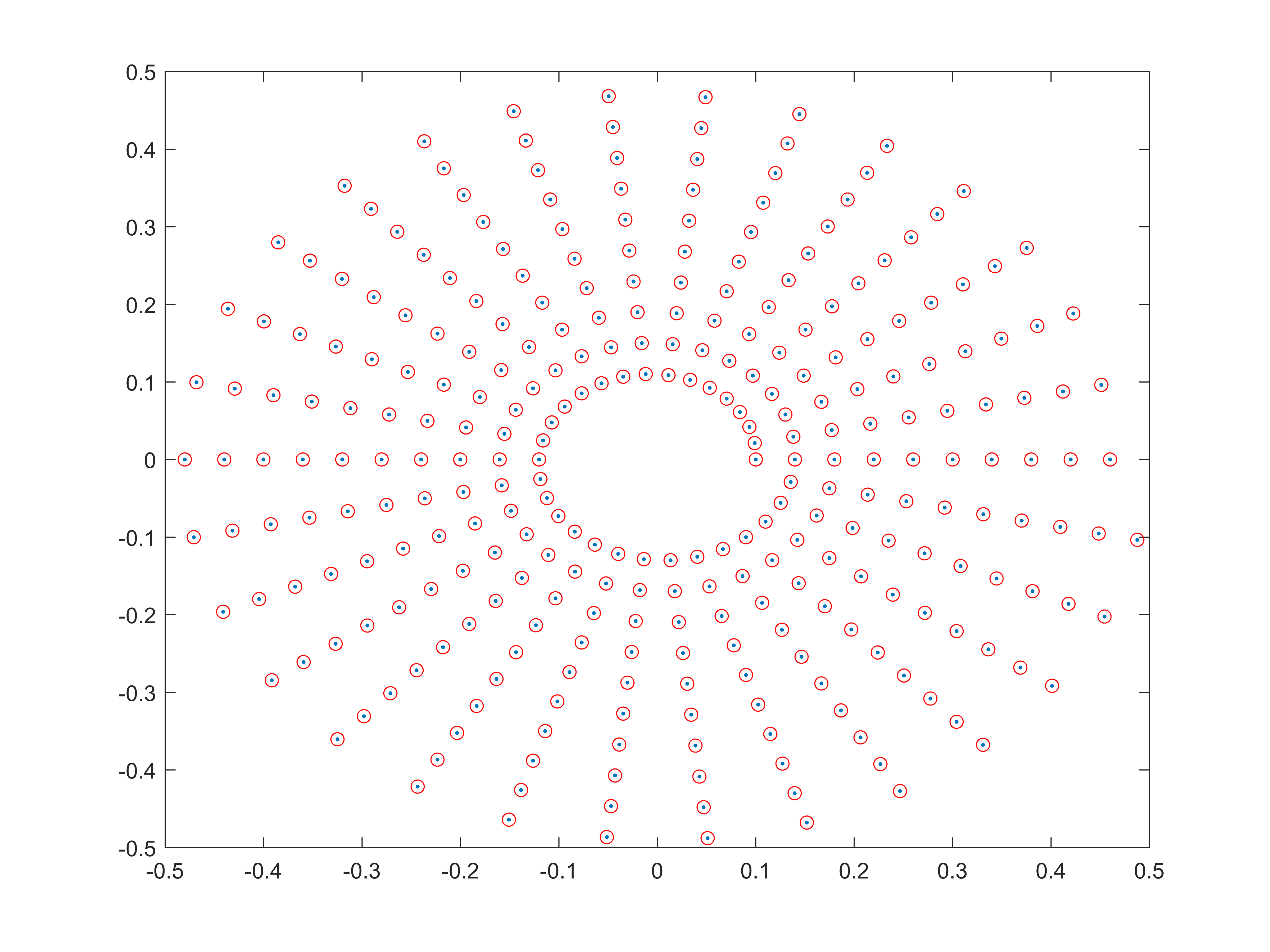}}
	\subfloat[][Error between estimated and true frequencies]{\includegraphics[width=0.47\linewidth]{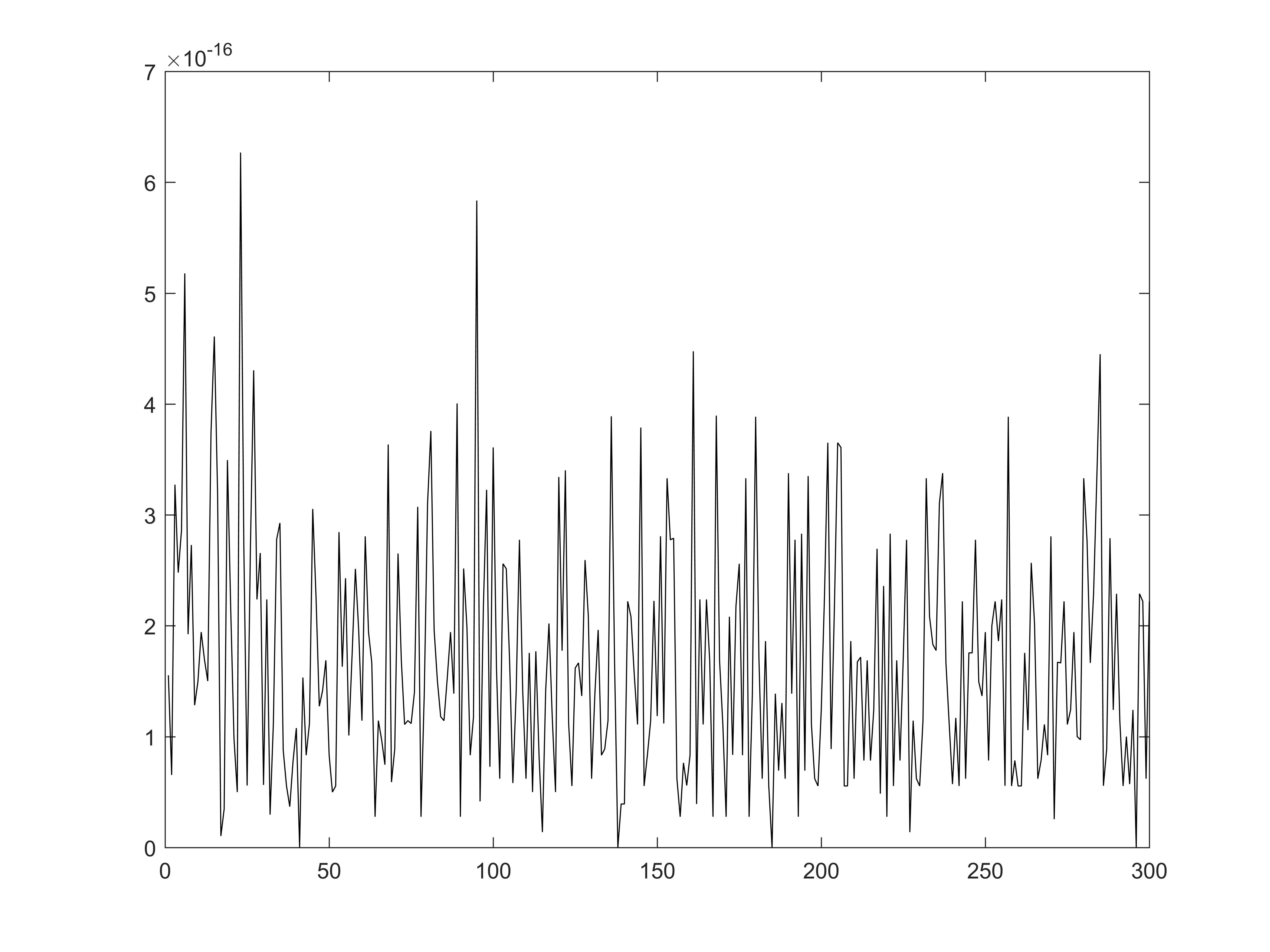}}
	\caption{\label{fig2dblock}Two dimensional example where a function $f$ of two variables is a linear combination of 300 exponential functions, which is sampled on a $61\times 61$-lattice (i.e. $N=31$).}
\end{figure}

There have been several suggestions on how to extend the one-dimensional results to several variables, c.f.~\cite{andersson2010nonlinear,cuyt1983multivariate, cuyt2011sparse,golyandina2013multivariate,harmouch2017structured,kung, kunis2016multivariate,peterreconstruction, rouquette2001estimation,moses,sauer2016prony,vanpoucke1994efficient,zhou2007novel}. See \cite{hua1992estimating} for an account of earlier attempts. The multidimensional case poses several difficulties; Let us assume that $f$ is of the form \begin{equation}
  \label{fmd}f(x)=\sum_{k=1}^K c_k e^{\zeta_k \cdot x},\quad c_k\in\C,~\zeta_k\in\C^d,
\end{equation}
(possibly distorted by noise) and that we wish to retrieve the frequencies $\zeta_k=(\zeta_{k,1},\ldots,\zeta_{k,d})$, where $\zeta_k \cdot x=\sum_{j=1}^d \zeta_{k,j}x_j.$ We suppose for the moment that $f$ has been sampled on an integer multi-cube $\{1,\ldots,2N-1\}^d$. A common engineering approach is to use some reformulation of the two dimensional problem to one-dimensional problems and from there recover the two dimensional components, see e.g.~the recent articles \cite{sward2016computationally,potts2013parametermult}. This causes a pairing problem that can lead to difficulties in practice, depending on the particular application at hand. An advantage with these type of methods is that they can be fast since they only sample data along lines.

To be more precise, one may average over all dimensions but one, say the $p$:th one, (or extract a ``fiber'' in dimension $p$ from the data) and then use a 1-d technique for estimating $\zeta_{1,p},\ldots,\zeta_{K,p}$, and repeat this for the remaining dimensions. However, first of all this limits the amount of frequencies we are theoretically able to retrieve to $N-1$ (in the noise free case, using e.g. ESPRIT, MUSIC or Prony's method), and secondly it is not clear how to pair the $Kd$ frequencies $\zeta_{k,p}$ into $K$ multi-frequencies $\zeta_k=(\zeta_{k,1},\ldots,\zeta_{k,d})$. Note that there are $K^d$ possible combinations, whereas only $K$ are sought. Also, these methods do not use the full data set and hence as an estimator they are likely to be more sensitive to noise, than a method which employs the full data set in the estimation. A related algorithm, which brings the idea of only using samples along certain lines to its extreme, is considered in \cite{plonka2013many} where it is shown that one can (theoretically) get away with $3K-1$ samples in total, in order to retrieve $K$ multi-frequencies (for $d=2$).

The focus in this article is different, we suppose that $f$ has been measured on a cube or more generally some irregular domain in $\R^d$, and propose a method which uses the full data set to retrieve the multi-frequencies and which avoids the pairing problem. This method is an extension of the ESPRIT-method. It was first introduced in \cite{rouquette2001estimation} in the 2d square case, and has recently been extended to several variables independently by \cite{sahnoun2017multidimensional} and \cite{steinwandt2016performance}. These articles also contain analysis of noise sensitivity and the former proposes accelerated versions. We review this method in its basic form in Sections \ref{secreview} (one variable) and \ref{block} (several variables).

The main contribution of the present work is to extend this method to work with data sampled on non-standard domains. This problem is also discussed in the 2d-setting in \cite{golyandina2013multivariate,Shlemov.Golyandina14conf-Shaped}. We also show how to realize block Hankel operators as multivariable summing operators, thereby providing a connection with theoretical results connecting the rank of these operators with functions of the form \eqref{fmd}, provided in \cite{andersson2016structure} and \cite{IEOT}.

\begin{figure}
	\centering
	\subfloat[][Real part of a 2d general domain data set constructed as a linear combination of 100 exponential functions.]{ \scalebox{1}[-1]{\includegraphics[width=0.47\linewidth,angle=-90]{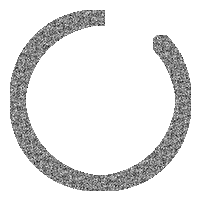} }}\\
	\subfloat[][True 2d frequencies in blue dots and estimated frequencies in red circles.]{	 \includegraphics[width=0.47\linewidth]{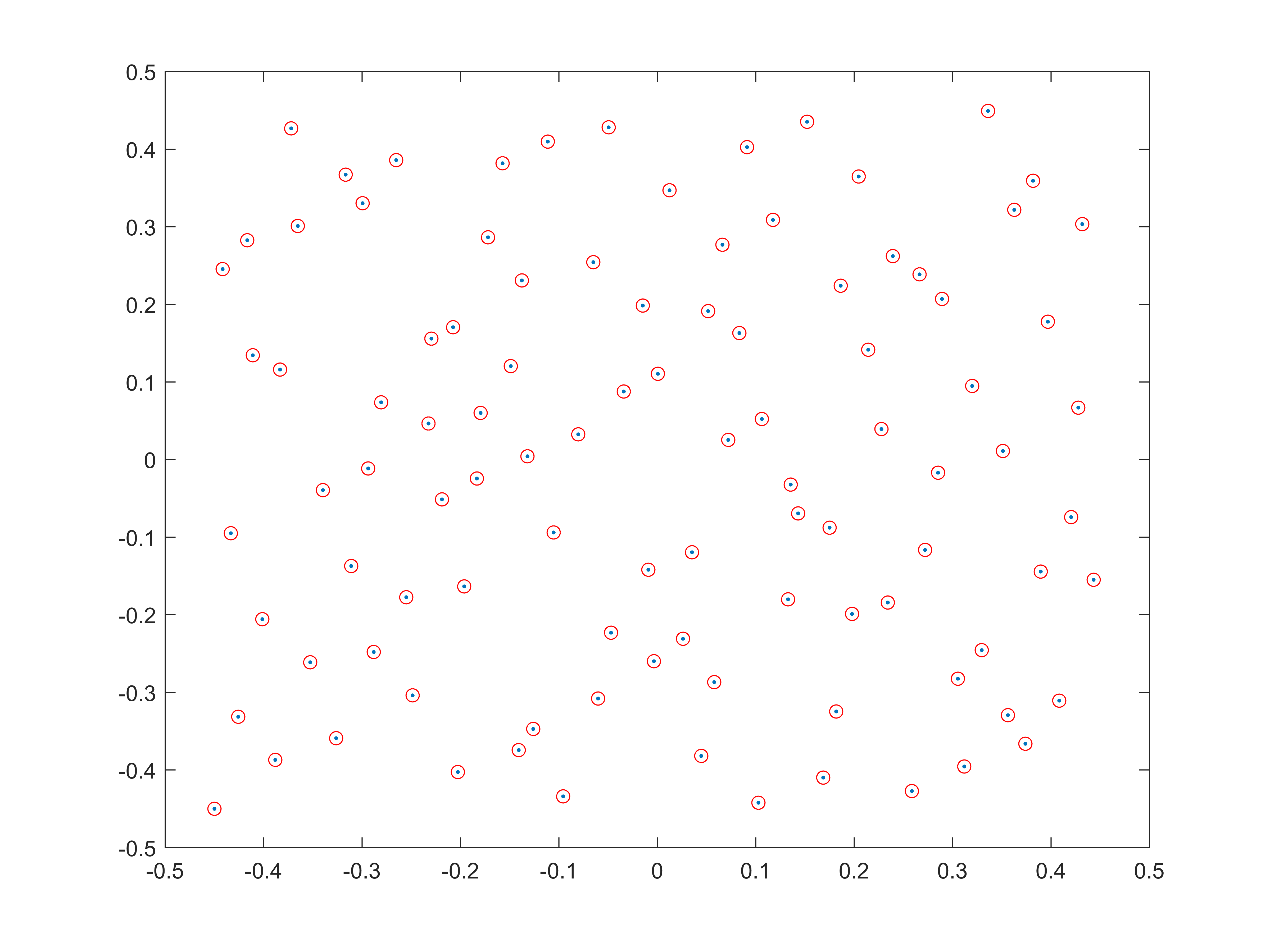}}
	\subfloat[][Error between estimated and true frequencies.]{	 \includegraphics[width=0.47\linewidth]{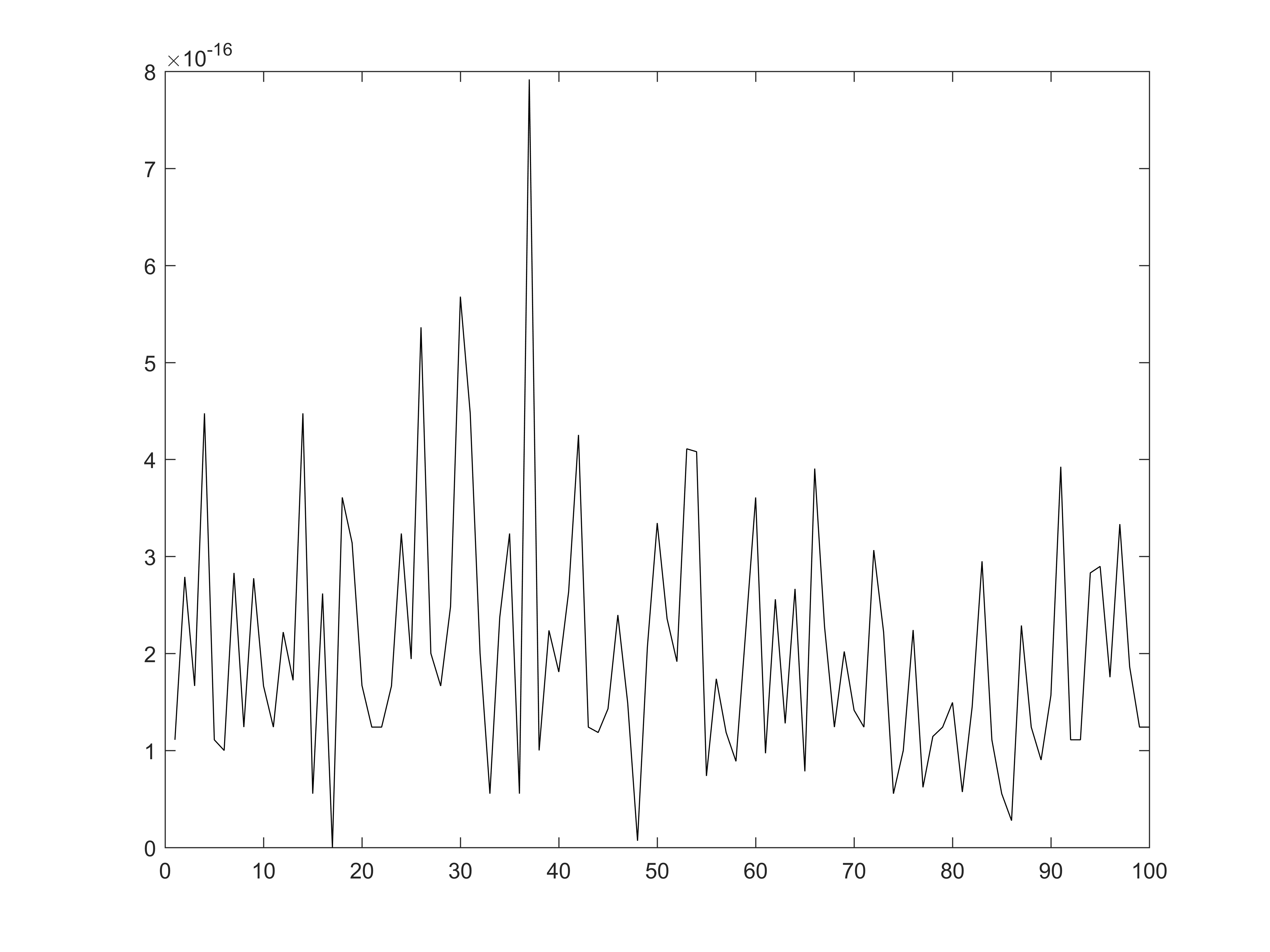}}
	\caption{\label{fig2dgen2}Two dimensional general domain example.}
\end{figure}

We now discuss a few examples to illustrate the pros and cons of the method suggested here. Suppose first that $f$ is sampled on a $19^3$ data-tensor, i.e. $N=10$ and $d=3.$ As it turns out, we are able to estimate $(N-1)\cdot N^{d-1}$ distinct multi-frequencies, avoiding degenerate cases. With the method considered here we may then retrieve 900 (randomly distributed) complex frequencies (assuming no noise), as opposed to 9 using reduction to one-dimensional techniques. By ``degenerate cases'' we for example rule out situations where many datapoints align parallel with one of the coordinate axes. For applications where this is the case the method considered here may not work. This is further discussed in Section \ref{block}, without making any rigorous attempt at formalizing what is meant by degenerate cases.

We now consider the example shown in Figure \ref{fig2dblock}, still in the multi-cube setting. The underlying function $f$ is in this case constructed as a linear combination of 300 purely oscillatory exponential functions with frequencies distributed along a spiral as depicted in panel c). The coefficients are chosen randomly, and $f$ is sampled on a 61 by 61 square grid (so $N=31$ and $d=2$). The real part of $f$ is shown in panel a). The block-Hankel matrix that is generated from $f$ is shown in panel  b). We have followed Algorithm \ref{alg_dd} to recover the underlying frequencies and they are recovered at machine precision error as shown in panel d).

Using reduction to one-dimensional techniques as discussed earlier, the maximum amount of frequencies that could have been retrieved from this data set is 30. With the present method it is $30*31=930$ (see \eqref{lok}) compared to the total amount of measurements $61^2=3721$. Or viewed from a different angle, to recover the 300 frequencies with reduction to one-dimensional techniques, we would theoretically need to have at least 601 measurements in each direction, and then it is not an easy problem to stably recover the 300 components in each dimension. Even if this is overcome, one then has $300\cdot 300=90000$ possible multi-frequencies to pick from, and 89700 has to be discarded, which may take some time. We refer again to \cite{sward2016computationally,potts2013parametermult} for two recent approaches on how to do this, while underlining that these methods may well be superior (and faster) in situations when fewer frequencies are sought and dense sampling along certain lines is not an obstacle in practice.

We now turn our focus to an example using the general domain Hankel structure, introduced in Section \ref{secdisc}, which allows us to work with samples of $f$ on domains with different shapes. In this case, we have used 100 exponential functions to generate the function $f$. This function is then sampled on a semi circular-like domain $\Omega$ shown in Figure \ref{fig2dgen2} a). Figure \ref{fig2dgen2} b) shows the distribution of the underlying frequencies, and Figure \ref{fig2dgen2} c) shows the reconstruction error in determining these frequencies. Again, the error is down at machine precision.

More examples are given in the numerical section, where we also test how the method performs with noise present, but we wish to underline that the main contribution of this paper is the deduction of an algorithm that is capable of exact retrieval in the noise-free case on general domains. Its performance as an estimator will be investigated elsewhere.

The paper is organized as follows, in Section \ref{secreview} we go over the essentials of the classical ESPRIT algorithm in one variable, in Section \ref{block} we extend this to the multi-cube case, using block Hankel matrices, and show how to solve the pairing problem. We then consider how to make this work for data sampled on non-cubical domains. This relies on so called \emph{general domain Hankel matrices}, which are introduced in Section \ref{secdi} and \ref{secdisc}. The subsequent extension of the ESPRIT-type algorithm from Section \ref{block} is given in Section \ref{fe5}. Numerical examples are given in Section \ref{numex}.

\section{One-dimensional ESPRIT}\label{secreview}

This section is intended as a quick review of ESPRIT in the simplest possible setting, to provide intuition for the more advanced versions introduced in later sections. For example we only work with square Hankel matrices (which we have found works best in practice), although non-square Hankel matrices will be considered in Section \ref{secdisc} and \ref{fe5} in the multivariable case. The study of ESPRIT and related algorithms in 1 dimension is of course a field in itself, and we refer e.g. to \cite{potts2013parameter} for a recent contribution with more flexibility.

Let us start by assuming that $f$ is a given one-dimensional function which is sampled at integer points, and that
$$
f(x)=\sum_{k=1}^K c_k e^{\zeta_k x},
$$
for distinct values of $\zeta_k$, $k=1, \dots K$.
A \emph{Hankel} matrix is a matrix $H$ which has constant values on its anti-diagonals. It is easy to see that each Hankel matrix can be generated from a function $h$ by setting $H(m,n)=h(m+n)$. Let us now consider the case where a (square $N\times N$) Hankel matrix $H$ is generated by the function $f$ above,  where $N>K$.

Let $\Lambda$ denote the Vandermonde matrix generated by the numbers $e^{\zeta_k}$, i.e., let
$$
\Lambda(j,k) = e^{\zeta_k j},\quad 1\leq j\leq N,~1\leq k\leq K
$$
and let $\Lambda_k$ denote the columns of $\Lambda$.
For the elements of the Hankel matrix $H$ it thus holds that
$$
H(m,n) = f(m+n) = \sum_{k=1}^K c_k e^{\zeta_k (m+n)} = \sum_{k=1}^K c_k e^{\zeta_k m} e^{\zeta_k n},
$$
implying that $H$ can be written as
$$
H = \sum_{k=1}^K c_k \Lambda_k \Lambda_k^T = \Lambda \mathrm{diag}(c) \Lambda^T,
$$
where $c=(c_1, \dots, c_K)$. From this observation it is clear that the rank of $H$ is $K$.

Now, by the singular value decomposition it follows that we can write $H=U\Sigma V^T$ where $U$ and $V$ are $N\times K$ matrices and $\Sigma$ a $K\times K$ diagonal matrix containing the non-zero singular values of $H$. It then holds that the corresponding singular vectors $U$ (or $V$) are linear combinations of the columns of $\Lambda$, i.e.,
$$
U =  \Lambda B,
$$
where $B$ is an invertible $K\times K$ matrix.

\begin{figure}
	\centering
\includegraphics[trim = 4cm 0cm 4cm 0cm,width=0.3\linewidth]{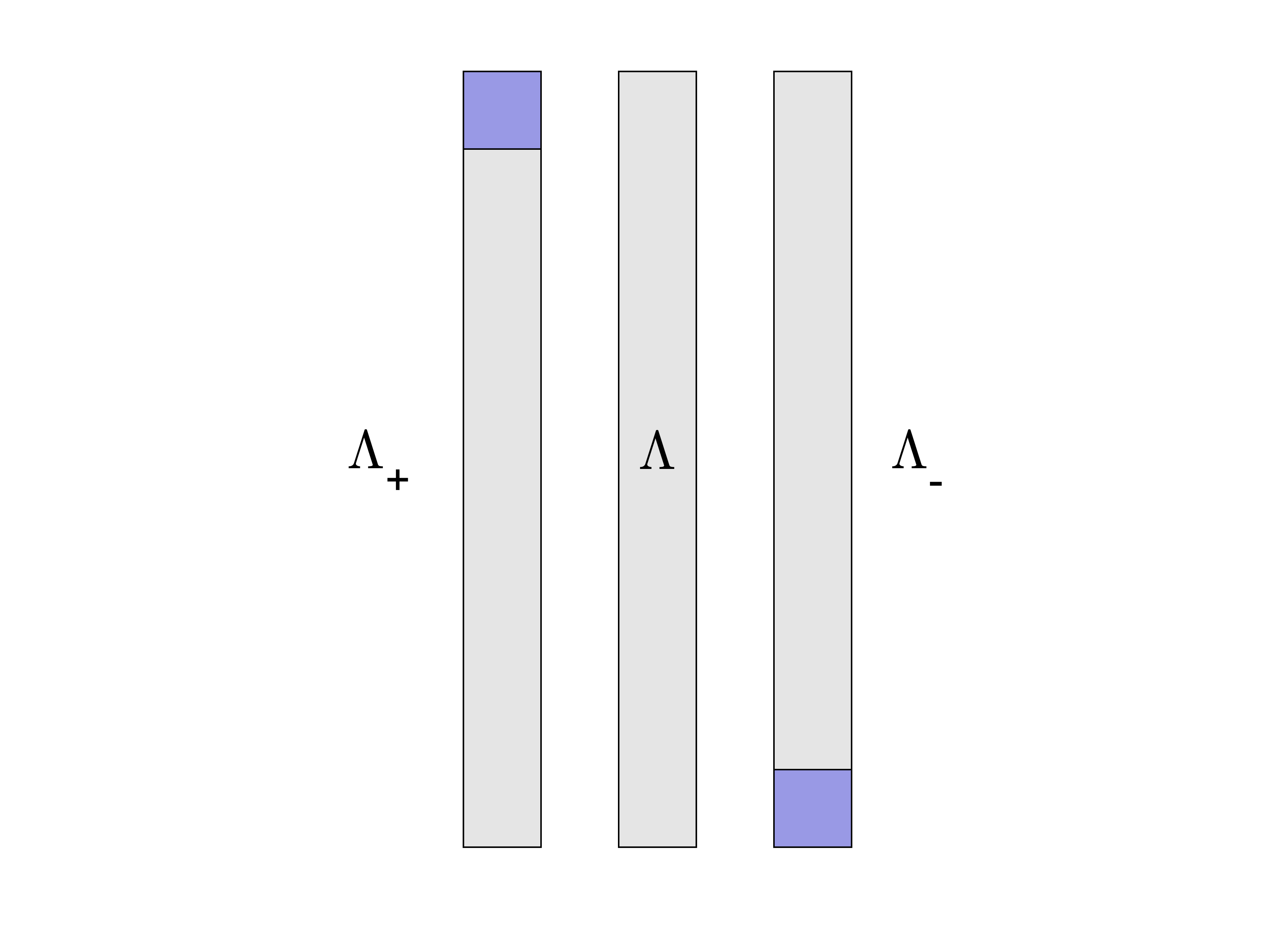}
  \caption{Illustration of the notation $\Lambda_+$ and $\Lambda_-$. For $\Lambda_+$ the first row is deleted, and for  $\Lambda_-$ the last row is deleted. The elements marked in blue are to be deleted.}
\end{figure}
Let $D = \mathrm{diag}(e^{\zeta_1},\dots, e^{\zeta_K})$. Let us also introduce the notation $\Lambda_+$ for the matrix that is obtained by deleting the first row from the matrix $\Lambda$, and similarly, let $\Lambda_-$ denote the matrix that is obtained by deleting the last row from $\Lambda$. There is a simple relationship between $\Lambda_+$  and $\Lambda_-$, namely
$$
\Lambda_+ =\Lambda_- D.
$$
Moreover, it holds that
\begin{align*}
&U_+ = \Lambda_+ B = \Lambda_- D B \\
&U_- = \Lambda_- B
\end{align*}
Since $U^*U=I$ and $N>K$, it typically holds that $U_-$ is injective. In particular it then has a left inverse $U_-^{\dagger}$ given by $$U_-^{\dagger} =   (U_{-}^\ast   U_{-})^{-1} U_{-}^\ast=(B^\ast \Lambda_{-}^\ast  \Lambda_{-} B)^{-1} B^\ast \Lambda_{-}^\ast
=B^{-1} (\Lambda_{-}^\ast  \Lambda_{-})^{-1} \Lambda_{-}^\ast .$$ Let us now consider
\begin{align*}
&A=U_-^{\dagger} U_+ =  B^{-1} (\Lambda_{-}^\ast  \Lambda_{-})^{-1} \Lambda_{-}^\ast \Lambda_- D B = B^{-1} D B.
\end{align*}
It follows that the columns of $B^{-1}$ are eigenvectors of the matrix $A$ with $e^{\zeta_k}$ being the corresponding eigenvalues.
This implies that the exponentials $e^{\zeta_k}$ can be recovered from $H$ by diagonalizing the matrix $A$, and this is the essence of the famous ESPRIT-method.

\begin{algorithm}
	\caption{One-dimensional ESPRIT}\label{alg_1d}
	\begin{algorithmic}[1]
		\State Form Hankel matrix $H(m,n)=f(m+n)$ from samples of $f$.
        \State Compute singular value decomposition $H =U \Sigma V^T$, retaining only non-zero singular values and corresponding vectors.
        \State Form $U_+$ and $U_-$ by deleting the first and last rows from $U$, respectively.
        \State Form $A = ( U_{-}^\ast   U_{-})^{-1} U_{-}^\ast U_{+}$.
        \State Diagonalize $A=B^{-1} \dia(\lambda_1, \dots, \lambda_K) B$ by making an eigenvalue decomposition of $A$.
        \State Recover $\zeta_k=\log(\lambda_k)$
	\end{algorithmic}
\end{algorithm}

We end by making a remark about the connection between Hankel matrices with low-rank and functions that are sums of exponential functions. For every function $f$ being a sum of $K$ (distinct) exponential functions the rank of the corresponding Hankel matrix is $K$ (given sufficiently many samples). The reverse state will typically be true, but there are exceptions. These exceptions are of degenerate character, and therefore this problem is often discarded in practice. For a longer discussion of these issues, see Section 2 and 11 in \cite{andersson2016structure}.

\section{Multi-dimensional ESPRIT for block-Hankel matrices}\label{block}

We now consider the case where $f$ is a given $d$-dimensional function which is sampled at integer points, and that
$$
f( x)=\sum_{k=1}^K c_k e^{{\zeta_{k}} \cdot {x}},
$$
where ${\zeta}_k=(\zeta_{k,1}, \dots, \zeta_{k,d})$ and $ x = (x_1, \dots, x_d)$. A $d$-dimensional Hankel operator can be viewed as a linear operator on the tensor product $\C^{N}\otimes\ldots\otimes\C^{N}$ whose entry corresponding to multi-indices $ m=(m_1, \dots, m_d)$ and $ n=(n_1, \dots, n_d)$ equals $f(m + n).$ This in particular entails that samples of $f$ needs to be available on a $(2N-1)^d$-dimensional grid. The restriction that we use the same amount of sample points in each direction is made for simplicity, the rectangular case is a special case of the more general setup in Section \ref{fe5}.

If we identify $\C^{N}\otimes\ldots\otimes\C^{N}$ with $\C^{N^d}$ using the reverse lexicographical order, i.e. by identifying entry $m$ in the former with \begin{equation}\label{o}m_1+m_2 N+m_3 N^2 + \dots m_d N^{d-1}\quad (0\leq m_j<N)\end{equation} in the latter, the corresponding operator can be realized as a block-Hankel matrix, which we denote by $\bb H$. Analogously, given multi-index $m$ we will write $\bb m$ for the number \eqref{o}, and if $u$ is an element in $\C^{N}\otimes\ldots\otimes\C^{N}$ we write $\bb u$ for its vectorized version in $\C^{N^d}$.


Just as for the one-dimensional case it can be shown that if
$\Lambda$ is the Vandermonde matrix generated by the vectors $\zeta_k$, i.e.,
$$
\bb{\Lambda}(\bb m,k) = e^{ \zeta_k \cdot m},
$$
then
$$
\bb H = \sum_{k=1}^K c_k \bb{\Lambda}_k \bb{\Lambda}_k^T = \bb{\Lambda} \mathrm{diag}(c) \bb{\Lambda}^T,
$$
and in particular $\bb H$ and $\bb \Lambda$ typically have rank $K$. If, due to insufficient sampling or a particular alignment of the $\zeta_k$'s, this affirmation is false, then the method described below will not apply, at least without further refinements. For example, maximal rank (i.e. $K$) can be achieved by rotation of the grid and/or by increasing $N$, which is further discussed in Section \ref{fe5}. Other tricks to deal with this issue is found in Section III of \cite{rouquette2001estimation}. Below we assume that problems related to rank deficiency or higher multiplicity of eigenvalues are absent, and discuss how to treat such issues towards the end of the section.

We thus assume that $\bb H$ and $\bb\Lambda$ has rank $K$ and write $\bb{H}=\bb{U}\Sigma \bb{V}^T$ for the singular value decomposition where we omit singular values that are zero and corresponding singular vectors, so in particular $\Sigma$ is a diagonal $K\times K$-matrix. As in the one-dimensional case it then holds that the singular vectors in $\bb{U}$ (or $\bb{V}$) are linear combinations of the columns of $\bb{\Lambda}$, i.e.,
\begin{equation}\label{ytr}
\bb{U} =  \bb{\Lambda} B,
\end{equation}
where $B$ is an invertible $K\times K$ matrix like before.

We now generalize the previous operation of deleting the first and last row respectively to several dimensions. By $\bb{\Lambda}_{p+}$ we denote the matrix that is obtained by deleting all entries of $\bb{\Lambda}$ with indices $\bb m$ corresponding to multi-indices $m$ of the form $(m_1, \dots 1,\dots m_d)$, i.e., all the first elements with respect to dimension $p$. Similarly, we denote by $\bb{\Lambda}_{p-}$ the matrix that is obtained by deleting from $\bb{\Lambda}$ all the elements related to multi-indices $(m_1, \dots N_p,\dots m_d)$, i.e., all the last elements with respect to dimension $p$. Let
\begin{equation*}
D_p=\dia{(e^{\zeta_{1,p}}, \dots\, e^{\zeta_{K,p}})}
\end{equation*}
and note that it in each dimension $p$ we have
$$
\bb{\Lambda}_{p+} =\bb{\Lambda}_{p-} D_p.
$$
\begin{figure}
	\centering
	\includegraphics[trim = 3cm 1cm 2cm 1cm,width=0.4\linewidth]{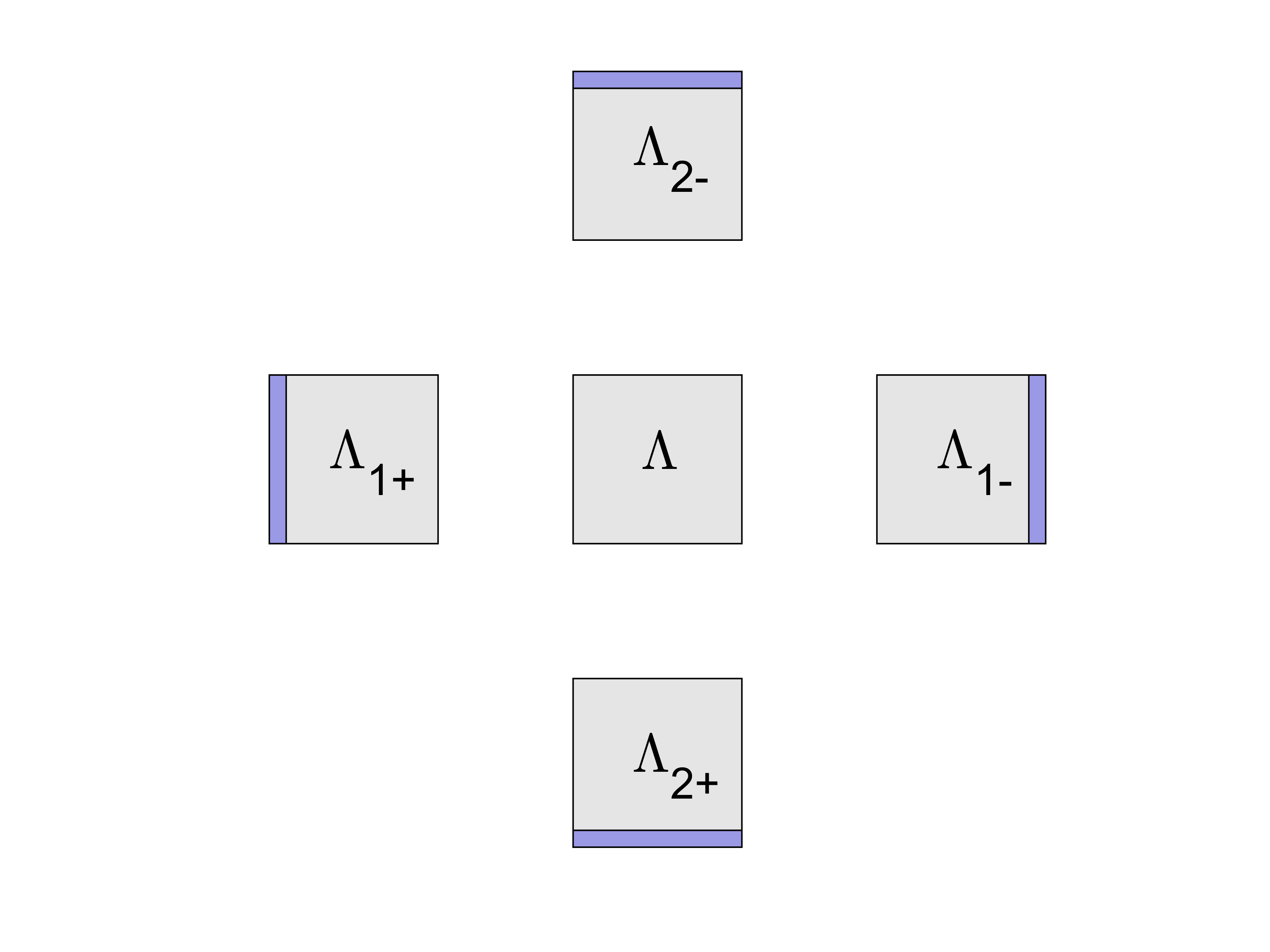}\hspace{2.5cm}
	\includegraphics[trim = 3.5cm 2cm 3.5cm 2cm,width=0.4\linewidth]{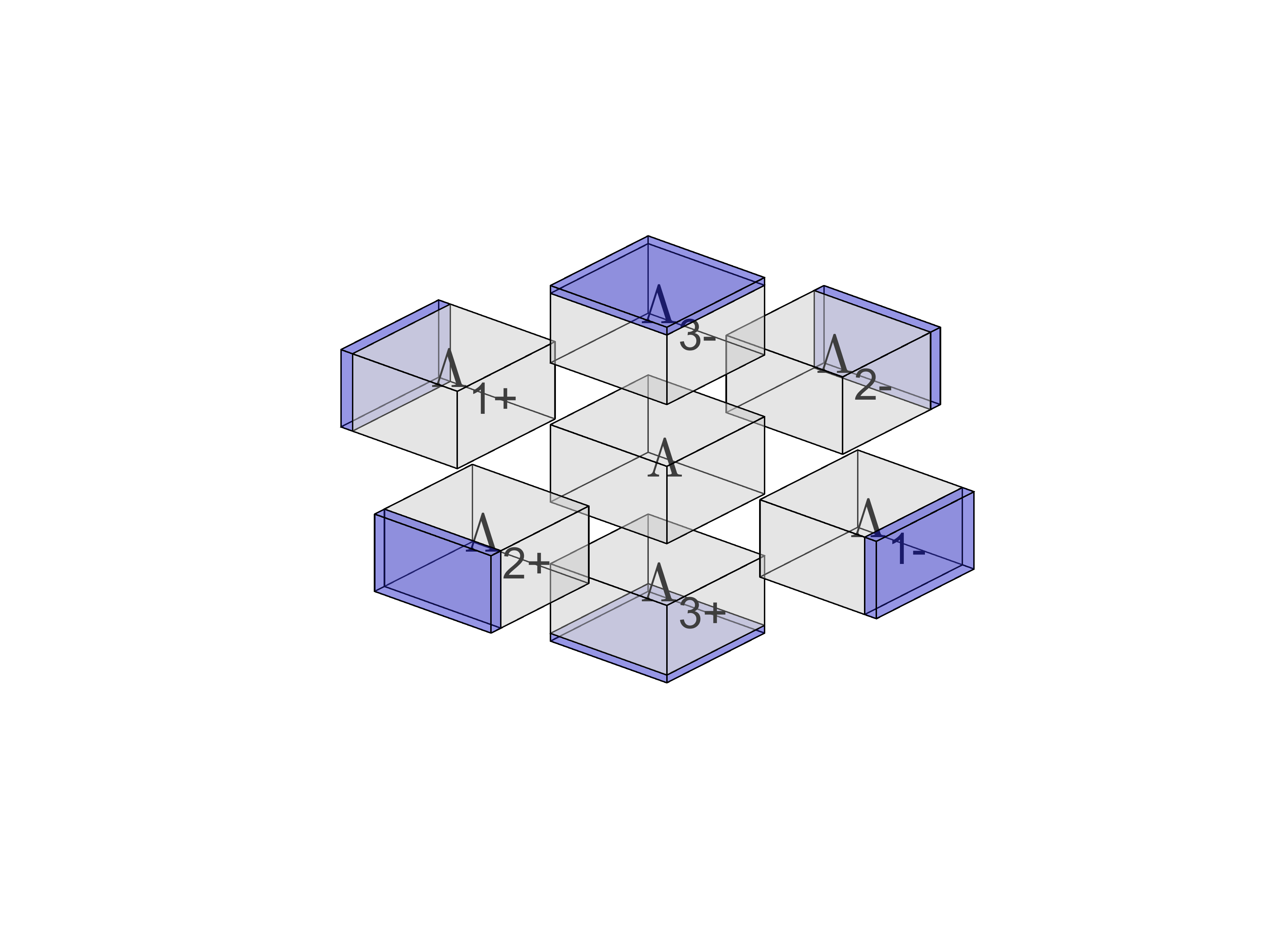}
	\caption{Illustration of the notation $\Lambda_{p+}$ and $\Lambda_{p+}$. For $\Lambda_{p+}$ the first elements in dimension $p$ are deleted, while for  $\Lambda_-$ the last elements in dimension $p$ are deleted. The elements marked in blue are to be deleted. The left panel shows a two-dimensional case, and in the right panel a three-dimensional case is shown.}
\end{figure}
Moreover, it holds that
\begin{align*}
&\bb{U}_{p+} =\bb{\Lambda}_{p+} B = \bb{\Lambda}_{p-} D_p B \\
&\bb{U}_{p-} = \bb{\Lambda}_{p-} B
\end{align*}
Let us now again consider
\begin{align*}
& A_p=\bb{U}_{p-}^\dagger \bb{U}_{p+} =  ( \bb{U}_{p-}^\ast   \bb{U}_{p-})^{-1} \bb{U}_{p-}^\ast   \bb{U}_{p+} \\&= (B^\ast \bb{\Lambda}_{p-}^\ast  \bb{\Lambda}_{p-} B)^{-1} B^\ast \bb{\Lambda}_{p-}^\ast  \bb{\Lambda}_{p-}  D_p B
\\&=B^{-1} (\bb{\Lambda}_{p-}^\ast  \bb{\Lambda}_{p-})^{-1} (B^{\ast})^{-1} B^{\ast} \bb{\Lambda}_{p-}^\ast  \bb{\Lambda}_{p-} D_p B \\&= B^{-1} D_p B.
\end{align*}
The remarkable observation here is that $B^{-1}$ \emph{simultaneously} diagonalizes all the matrices $A_p$, and that the eigenvalues to $A_p$ are $(e^{\zeta_{1,p}}, \dots\, e^{\zeta_{K,p}})$ in the correct order independent of $p$. This implies that all the exponentials $e^{\zeta_{k,1}}$, $k=1,\ldots ,K$, can be recovered from $H_f$ by diagonalization of $A_1$, and then we may use \textit{the same} eigenvectors for diagonalization of the remaining $A_p$'s, and in this way obtain the multi-frequencies $ \zeta_k=(\zeta_{k,1},\ldots,\zeta_{k,d})$ for $k=1,\ldots ,K$ directly without any need for a pairing procedure. This phenomenon was first observed in \cite{rouquette2001estimation} in the two-dimensional case, although it is not very clear to see. A more accessible account, which extends \cite{rouquette2001estimation} to the higher dimensional case, is given in \cite{sahnoun2017multidimensional} and independently in \cite{steinwandt2016performance}. To summarize, the algorithm in its most simple form reads as follows:

\begin{algorithm}
	\caption{$d$-dimensional ESPRIT}\label{alg_dd}
	\begin{algorithmic}[1]
		\State Form the $d$-block Hankel matrix $\bb{H}(\bb m,\bb n)=f( m+ n)$ from samples of $f$ .
		\State Compute the singular value decomposition $\bb{H} =\bb{U} \Sigma \bb{V}^T$.
		   \For {$p=1, \dots, d$}
		   \State Form $\bb{U}_{p+}$ and $\bb{U}_{p-}$ by deleting first and last elements in dimension $p$, respectively.
		   \State Form $A_p =  ( \bb{U}_{p-}^\ast  \bb{U}_{p-})^{-1} \bb{U}_{p-}^\ast   \bb{U}_{p+}$.
		   \EndFor		
		\State Diagonalize $A_1=B^{-1} D_1 B$ by making an eigenvalue decomposition of $A_1$.
	    \For {$p=1, \dots, d$}
		   \State Compute the diagonal matrices $D_p = \dia(\lambda_{1,p}, \dots, \lambda_{K,p})=B A_p B^{-1}$.
		   \State Recover $\zeta_{k,p}=\log(\lambda_{k,p})$ (which are automatically correctly paired).
		\EndFor		
	\end{algorithmic}
\end{algorithm}

We now discuss problems that may arise in the above approach. The algorithm works as stated as long as $(\bb{\Lambda}_{p-})^\ast   \bb{\Lambda}_{p-}$ is invertible for all $p$, and moreover we need that no eigenvalue in $D_1$ has multiplicity higher than 1 in order for the eigenvectors in $B$ to be correctly determined. The first limitation leads to the restriction \begin{equation}\label{lok}K\leq N^{d-1}(N-1),\end{equation} since $\bb{\Lambda}_{p-}$ needs to have fever columns than rows. This restriction was mentioned already in the introduction and also appears in Lemma 2 of \cite{sahnoun2017multidimensional}, where it is noted that this condition holds generically if the frequencies are sampled at random. In the same article, a linear combination step is also applied to avoid problems with multiplicity. Indeed, suppose for simplicity that $d=2$ and that $D_1$ does not have distinct eigenvalues. Then \begin{equation}\label{p0}\alpha A_1+\beta A_2=B^{-1}(\alpha D_1+\beta D_2)B\end{equation} has distinct eigenvalues for most choices of $\alpha$ and $\beta$, so the problem can be circumvented by randomly choosing $\alpha,\beta$ and compute $B$ from the above linear combination. This trick is also employed in \cite{potts2013parametermult}, albeit for a different algorithm.

A final remark on Algorithm \ref{alg_dd} concerning time-complexity. In \cite{sahnoun2017multidimensional} it is noted that Algorithm \ref{alg_dd} can be slow since performing an SVD on a large matrix is time consuming. It is suggested to circumvent this problem by using ``truncated SVD'' which computes the $K$ first singular vectors based on a variant of Lanczos algorithm. We note that similar improvements can be applied to the Algorithm \ref{alg_gendomain} presented below, but we do not follow these threads here.

This completes the multidimensional version of ESPRIT in the case when data $f$ is sampled on a regular multi-cube. Next we address the more general setting of data measured on regular grids with various shapes, which is the main contribution of this paper, presented in Algorithm \ref{alg_gendomain}. For this we need to introduce so called general domain Hankel matrices, by viewing block Hankel matrices as multi-dimensional summing operators.

\section{The summation operator formalism in one variable}\label{secdi}
\begin{figure} \centering
	\includegraphics[width=0.5\textwidth,height=0.37\textwidth,trim={0cm 0cm 0cm 0cm}]{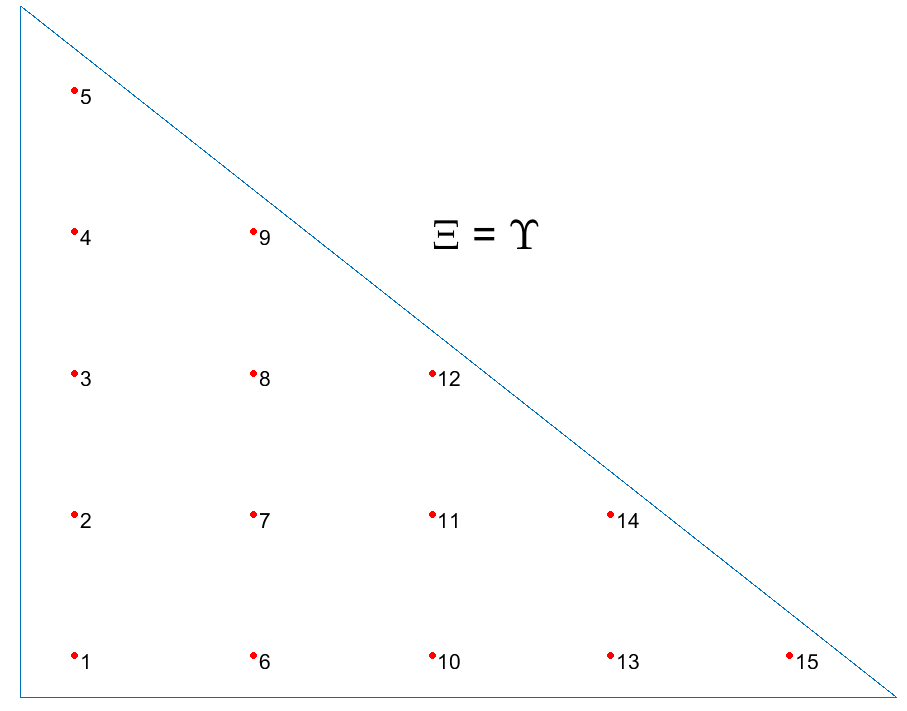}\caption{\label{fig:domainXi} Samples of a triangular domain with corresponding ordering. If the vertical axis is considered to be the first coordinate, then this ordering is the reverse lexicographical ordering.}
\end{figure}

Given a subset $\Upsilon\subset \N$, let $\ell^2(\Upsilon)$ be the ``sequences'' indexed by $\Upsilon$ and equipped with the standard $\ell^2-$norm. If $\Upsilon=\{1,2,\ldots,N\}$, then $\ell^2(\Upsilon)$ reduces to the standard $\C^N$. Note that we may define a classical Hankel matrix $H$ given by the sequence $f=(f_2,f_3,\ldots,f_{2N})$ as the operator with domain and codomain equal to $\ell^2(\{1,\ldots,N\})$, given by the summation formula
\begin{equation}\label{sumop1}
(H_{{f}}({a}))_n=\sum_{m\in \{1,\ldots,N\}} f_{m+n}a_m,\quad n\in \{1,\ldots,N\}
\end{equation}
where ${a}$ represents any element in $\ell^2(\{1,\ldots,N\})$. This suggests the following generalization; Let $M,N\in\N$ such that $M+N=2P$ be given, and consider $H_{{f},M,N}:\ell^2(\{1,\ldots,M\})\rightarrow\ell^2(\{1,\ldots,N\})$ given by
\begin{equation}\label{sumop2}
(H_{{f},M,N}({a}))_n=\sum_{m\in \{1,\ldots,M\}} f_{m+n}a_m,\quad n\in \{1,\ldots,N\}.
\end{equation}
We remark that \begin{equation}\label{domain1}
\{1,\ldots,M\}+\{1,\ldots,N\}= \{2,\ldots,2P\},
\end{equation}
where the latter is the grid on which $f$ is sampled. A moments thought reveals that $H_{{f},M,N}$ equals the $N\times M$ Hankel matrix given by ${f}$. For example, if ${f}=(1,2,3,4,5)$, then $P=3$ and we can pick e.g. $M=4$ and $N=2$. We then have
\begin{equation}\label{hankelex}H_{{f},4,2}=\left(
                                                                                 \begin{array}{ccccc}
                                                                                   1 & 2 & 3 & 4 \\
                                                                                    2& 3 & 4 & 5  \\                                                                                   \end{array}
                                                                               \right).
\end{equation}
As long as $M$ and $N$ are larger than $K$, it is possible to perform ESPRIT on these rectangular Hankel matrices as well, although square ones seems to be preferred.

\section{General domain Hankel operators}\label{secdisc}

\begin{figure}
\centering 	
{\includegraphics[width=0.48\textwidth]{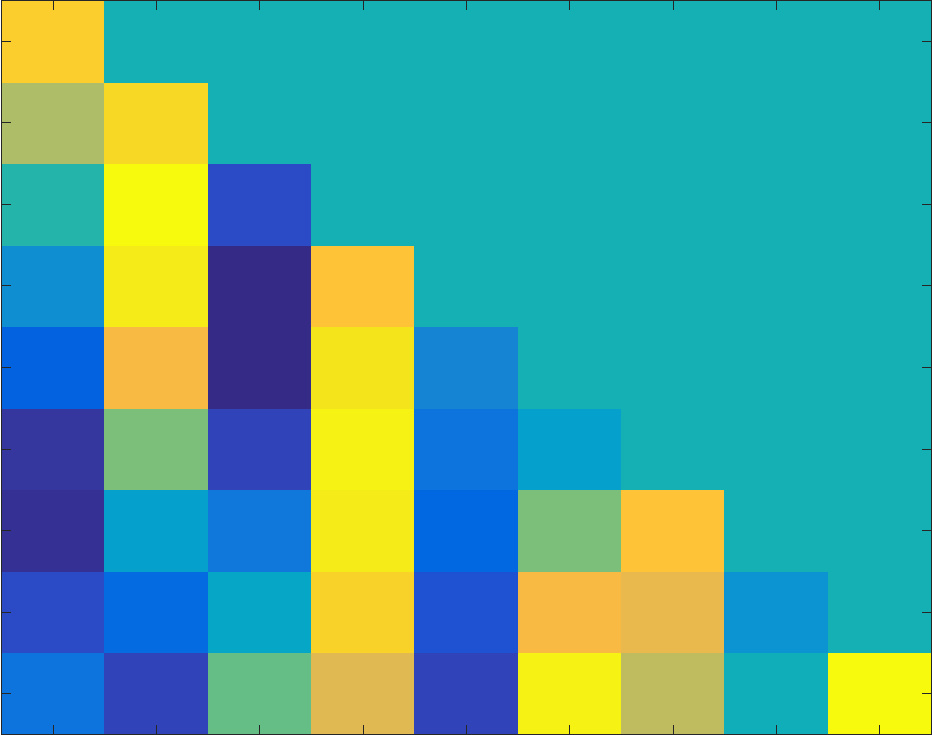}}
{\includegraphics[width=0.48\textwidth]{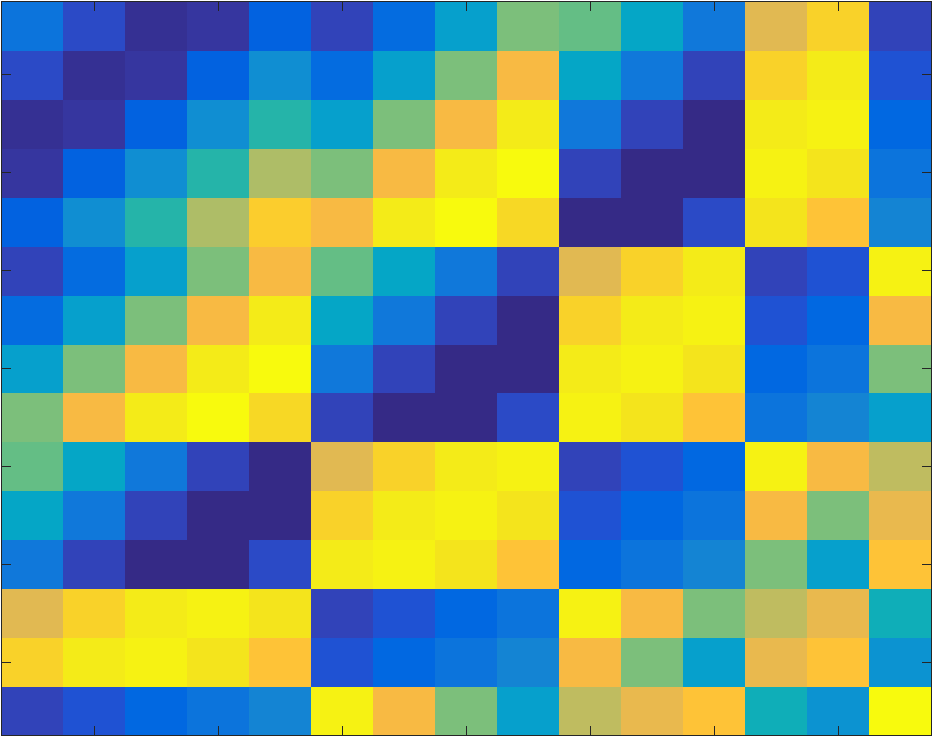}}
\caption{\label{fig:conv_support} a) Function defined on $\Omega=\Xi+\Upsilon$ with $\Xi=\Upsilon$ as in Figure \ref{fig:domainXi}. b) The corresponding general domain Hankel matrix, using the ordering from Figure \ref{fig:domainXi}. }
\end{figure}

General domain Hankel integral operators were introduced in \cite{IEOT} (albeit under the name truncated correlation operators) and their discretizations were studied in \cite{andersson2016structure} where the term \emph{general domain Hankel matrix} was coined. The discrete versions have been studied earlier in e.g. \cite{mourrain2000multivariate} under the name quasi-Hankel matrices. We briefly revisit their construction, which is explained in greater detail in Section 4 of \cite{andersson2016structure}.

\begin{figure}
	\centering
	\setlength{\unitlength}{\linewidth}
	\subfloat[][General domain setup with $\Upsilon$ on top, $\Xi$ in the middle and $\Omega$ in the bottom. The $\Xi$ set is a $11\times 11$ grid]{\includegraphics[height=0.7\linewidth, width=0.45\linewidth,trim={7cm 0 7cm 0}]{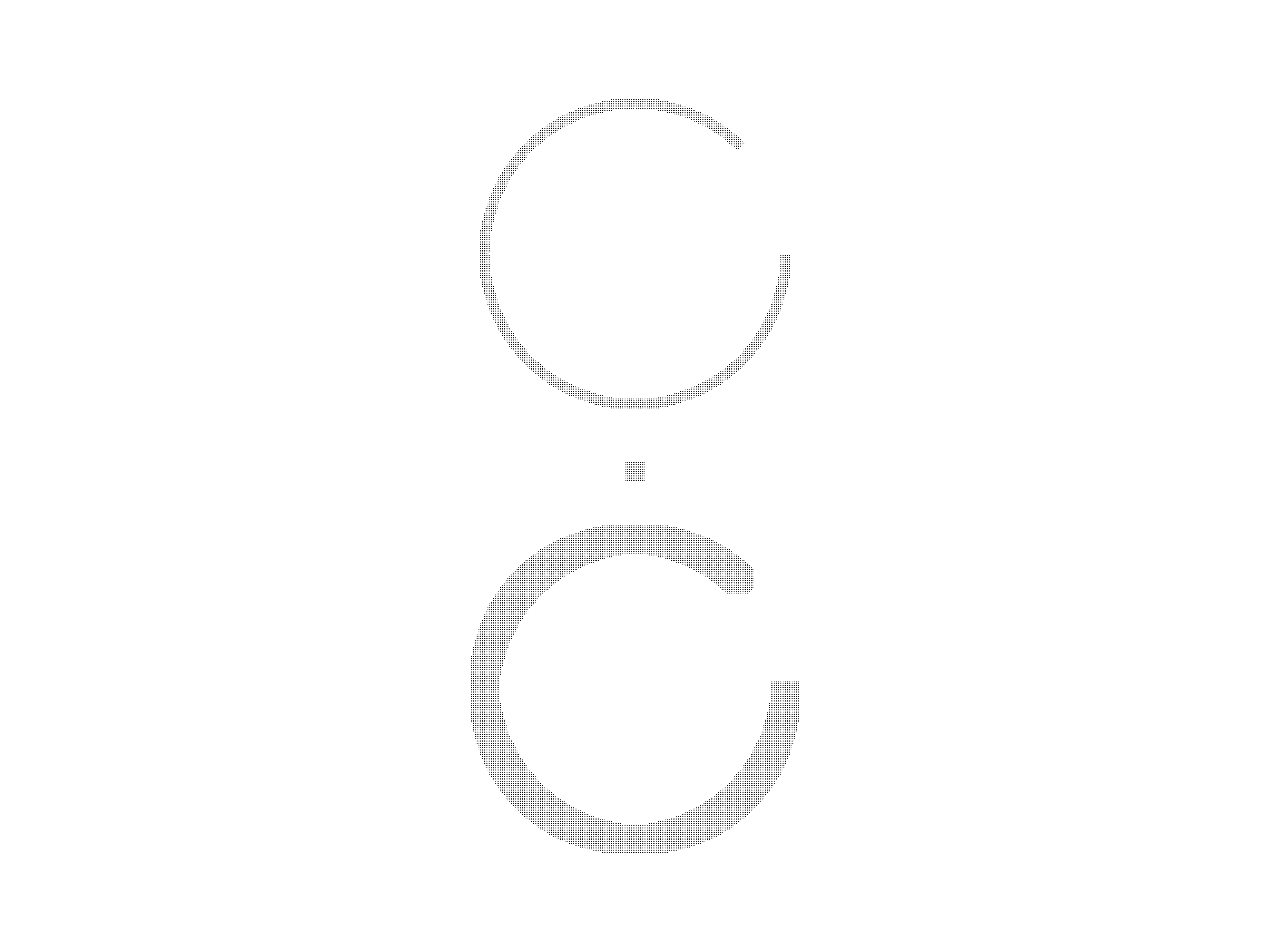}}
	\subfloat[][The corresponding general domain Hankel matrix]{\includegraphics[height=0.7\linewidth, width=0.45\linewidth]{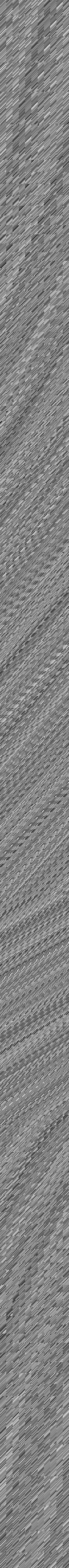}	}
	\caption{\label{fig2dgen1} Two dimensional general domain example.}
\end{figure}

Let $\Xi,\Upsilon$ be any bounded subsets of $\Z^d$, and in analogy with \eqref{domain1} set
\begin{equation}\label{domain2}
\Upsilon+\Xi= \Omega,
\end{equation}
as in Figure \ref{fig:domainXi} and \ref{fig2dgen1} a). Suppose we are interested in a function $f$ on some domain in $\R^d$, and that $\Upsilon$ and $\Xi$ are chosen such that $\{j_1\Delta x_1 ,\ldots j_d\Delta x_d): {j}\in\Omega\}$ cover the domain where $f$ is defined, for some choice of sampling length $\Delta x_1,\ldots,\Delta x_d$. The function $f$ thus gives rise to a \emph{multidimensional sequence}
 \begin{equation}\label{sample}{f}_{{j}}=f(j_1\Delta x_1 ,\ldots, j_d\Delta x_d),\quad {j}\in{\Omega},\end{equation} or more formally a function in $\ell^2(\Omega)$, see Figure \ref{fig:conv_support} a). We will refer to such functions as \emph{md-sequences}, to distinguish them from ordinary sequences, (i.e. vectors in $\C^N$). In analogy with \eqref{sumop2}, the md-sequence $f$ gives rise to a corresponding general domain Hankel operator \begin{equation}\label{sumop} (H_{{f},\Upsilon,\Xi}(g))_{n}=\sum_{m\in {\Upsilon}}f_{m+n}g_{m},\quad n\in {\Xi},\end{equation} where $g$ is an arbitrary md-sequence on ${\Upsilon}$. When $\Upsilon,\Xi$ are clear from the context or irrelevant, we drop them from the notation.

We may of course represent $g$ as a vector by ordering the entries (as in Figure \ref{fig:domainXi} for a triangle). More precisely, by picking any bijection \begin{equation}\label{bij}o_y:\{1,\ldots,|{\Upsilon}|\}\rightarrow{\Upsilon}\end{equation}
(where $|\Upsilon|$ denotes the amount of elements in the set $\Upsilon$), we can identify $g$ with the vector $\bold{g}$ given by $$(\bold{g}_j)_{j=1}^{|{\Upsilon}|}=g(o_y(j)).$$
Letting $o_x$ be an analogous bijection for ${\Xi}$, it is clear that $H_{{f},\Upsilon,\Xi}$ can be represented
as a $|\Xi|\times |\Upsilon|$-matrix, where the $(n,m)$'th element is $f(o_x(n)+o_y(m))$, see Figure \ref{fig:conv_support} b) and \ref{fig2dgen1} b). Such matrices will be called \emph{general domain Hankel matrices} and denoted $\bold{H}_{f,\Upsilon,\Xi}$, letting the bijections $o_x$ and $o_y$ be implicit in the notation (usually we will use the reverse lexicographical order as in Section \ref{block}). In particular, if we set $$\Upsilon=\Xi=\{1,\ldots, N\}^d,$$ we retrieve in this way the block Hankel matrices discussed in Section \ref{block}. Note that in the one dimensional case, $H_{f,\Upsilon,\Xi}$ and $\bold{H}_{f,\Upsilon,\Xi}$ are virtually the same thing, whereas this is not the case in several variables. The former is an operator acting on md-sequences in $\ell^2(\Upsilon)$ to md-sequences in $\ell^2(\Xi)$, the latter is a matrix representation of the former that can be used e.g. for computer implementations.

An example where $\Xi=\Upsilon$ (and hence also $\Omega$) are triangles is shown in Figures \ref{fig:domainXi} and \ref{fig:conv_support}. The domains $\Xi=\Upsilon$ are shown in Figures \ref{fig:domainXi} with a particular ordering. $\Omega$ is then a triangle of side-length 9, and a function $f$ on $\Omega$ is shown in Figure \ref{fig:conv_support} a). The corresponding general domain Hankel matrix, given the ordering from Figure \ref{fig:domainXi}, is shown in b). Note especially how columns 3, 4 and 5 are visible and show up as rectangular Hankel matrices of different sizes. Columns 3 and 5 both give rise to a (blue) square Hankel matrix on the ``block-structure diagonal'', which also show up as rectangular Hankel matrices off the ``block-diagonal''. Column 4 generates the (yellow) Hankel matrices that are never on the diagonal, corresponding to the fact that every second number in a 1d Hankel matrix does not show up on the diagonal.

Figure \ref{fig2dgen1} a) shows how a semi-circular domain $\Omega$ is constructed from the two domains $\Upsilon$ and $\Xi$. Note that the set $\Xi$ is a rather small domain - a grid of size $11\times 11$. The size of $\Xi$ will determine how many frequencies that can recovered, in this case 11*10=110, (which follows by the same argument that led to formula \eqref{lok}). The general domain Hankel matrix that is constructed from $f$ is shown in Figure \ref{fig2dgen1} b).

\section{Frequency retrieval; general domain ESPRIT}\label{fe5}

%
%

As earlier, we assume that we have samples of a function $f$ of the form \begin{equation}\label{lokf}
f( x)=\sum_{k=1}^K c_k e^{{\zeta_{k}} \cdot {x}},
\end{equation}
where the axes have been scaled so that we sample at integer points. However, we now assume that samples are available only on $\Omega$ which is a domain of the form $\Xi+\Upsilon$ for some domains $\Xi$ and $\Upsilon$, where we only require that $\Xi$ is convex (in the sense that it arise as the discretization of a convex domain). Given $\lambda\in\C^{d}$ we let $\Lambda_{\Omega}(\lambda)$ denote the md-sequence $\lambda^j$ for $j\in\Omega$, where we use multi-index notation $\lambda^j=\lambda_{1}^{j_1}\lambda_{2}^{j_2}\ldots \lambda_{d}^{j_d}$. Rephrased, we suppose that $f$ can be written
\begin{equation}
  \label{expsepmddiscrete}{f}=\sum_{k=1}^K c_k {\Lambda}_{\Omega}(\lambda_k),\quad \lambda_k=(e^{{\zeta_{k,1}} },\ldots,e^{{\zeta_{k,d}} }).
\end{equation}
It is then easy to see that \begin{equation}\label{sunday}H_f(g)=\sum_{k=1}^K c_k H_{\Lambda_{\Omega}(\lambda_k)}(g)=\sum_{k=1}^K c_k \Lambda_{\Xi}(\lambda_k)\langle g,\overline{\Lambda_{\Upsilon}(\lambda_k)}\rangle,\end{equation}
which in particular shows that $H_f$ is a rank $K$ operator, assuming of course that the $c_k$'s are non-zero, the $\lambda_k$'s are distinct and that $\{\Lambda_{\Xi}(\lambda_k)\}_{k=1}^K$ and $\{\Lambda_{\Upsilon}(\lambda_k)\}_{k=1}^K$ form linearly independent sets.

A formal investigation of when this condition holds becomes very involved, but is possible given certain assumptions on $\Xi$ and $\Upsilon$, we refer to the work of B. Mourrain et al.~\cite{mourrain2000multivariate,brachat2010symmetric}. Another viewpoint is to consider general domain Hankel operators as discretizations of corresponding integral operators, whose rank structure is easier to characterize, we refer to the articles \cite{IEOT,andersson2016structure} (by the authors) as well as the recent contribution \cite{mourrain2016polynomial} by B. Mourrain. In particular, it is well known that any set of exponential functions with different exponents always is linearly independent, as functions on whatever fixed open set in $\R^n$. The condition that $\{\Lambda_{\Upsilon}(\lambda_k)\}_{k=1}^K$ be linearly independent will thus always be satisfied given that the sampling is dense enough. We refer to Section 6 of \cite{andersson2016structure} for more information on passing between discrete (summing) and continuous (integral) general domain Hankel operators. In particular, technical conditions (see (6.3)) on the boundary of the domains in $\R^d$ which correspond to $\Upsilon$ and $\Xi$ are given so that the sampled md-vectors converge, in a sense made precise, to the corresponding exponential functions (see the proof of Theorem 6.4 in \cite{andersson2016structure}). In the remainder of this paper, we assume that the linear independence condition is fulfilled.

Letting $o_x$ be the reverse lexicographical ordering on $\Xi$, we identify md-sequences $u$ in $\ell^2(\Xi)$ with sequences (vectors) $\bold{u}$ in $\ell^2(|\Xi|)$, as explained in Section \ref{secdisc}. The SVD of $\bold{H}_f$ thus gives rise to singular vectors $\bold{u}_k$, $\bold{v}_k$, whose md-sequence counterparts satisfy $$H_f(v_k)=\sigma_k u_k,\quad H_f^*(u_k)=\sigma_k v_k.$$ Since $\{u_k\}_{k=1}^K$ span the range of $H_f$, i.e.~the span of $\{\Lambda_{\Xi}(\lambda_k)\}_{k=1}^K$ by \eqref{sunday}, it follows that we can write \begin{equation}\label{dr}u_j=\sum_{k=1}^K b_{j,k}\Lambda_{\Xi}(\lambda_k)\end{equation} where the numbers $b_{j,k}$ form a square $K\times K$ invertible matrix. Let $\bold{U}$ be the matrix with columns $\bold{u}_1,\ldots,\bold{u}_K$ and let $\bold{\Lambda}$ be the matrix with the columns $\bold{\Lambda}_{\Xi}(\lambda_1),\ldots,\bold{\Lambda}_{\Xi}(\lambda_K)$. The relation \eqref{dr} can then be expressed \begin{equation}\label{dr1}\bold{U}= \bold{\Lambda} B\end{equation}
(compare with \eqref{ytr}).

By a ``fiber'' in $\Xi$ we refer to a subset obtained by freezing all variables but one. Since we have assumed that the grid $\Xi$ is the discretization of a convex domain, we have that each fiber in the first dimension is of the form $\{M_1,\ldots, M_2\}\times \{m_2\}\times \ldots\times\{m_d\}$, where $M_1\leq M_2$ depend on the ``frozen'' variables $m_2,\ldots,m_d$. We denote the grid that arises by removing the first (respectively last) element of each such fiber by $\Xi_{+,1}$ (respectively $\Xi_{-,1}$), where we assume that the sampling has been done so that no fiber consists of a singleton. Analogous definitions/assumptions apply to the other variables, yielding grids $\Xi_{\pm, 2},\ldots,\Xi_{\pm, d}$.

Now, given a fixed dimension $p$ and a md-sequence $w$ in $\ell^2(\Xi)$, we let $w_{\pm, p}$ be the md-sequence restricted to the grid $\Xi_{\pm, p}$. Moreover, given a matrix like $\bold{U}$, whose columns are given by the vectorizations of the  md-sequences $u_k$ in $\ell^2(\Xi)$, we denote by $\bold{U}_{+,p},$ (resp.~$\bold{U}_{-,p}$) the matrix formed by the vectorized md-sequences $u_{k,+,p}\in \ell^2(\Xi_{+,p})$ (resp. $u_{k,-,p}\in \ell^2(\Xi_{-,p})$). Equation \eqref{dr1} then implies that \begin{equation}\label{dr2}\bold{U}_{\pm,p}= \bold{\Lambda}_{\pm,p} B.\end{equation}
Setting $D_p=\dia{(e^{\zeta_{1,p}}, \dots\, e^{\zeta_{K,p}})}$, we also have
$\bold{\Lambda}_{+,p}=\bold{\Lambda}_{-,p}D_p$
which combined implies that \begin{equation}\label{dr5}\bold{U}_{+,p}= \bold{\Lambda}_{-,p}D_pB.\end{equation}
In analogy with the computations in Section \ref{block} we have
\begin{equation*}\label{gt}A_p:=(\bold{U}_{-,p}^*\bold{U}_{-,p})^{-1}\bold{U}_{-,p}^*\bold{U}_{+,p}=(B^*\bold{\Lambda}_{-,p}^*\bold{\Lambda}_{-,p}B)^{-1}B^*
\bold{\Lambda}_{-,p}^*\bold{\Lambda}_{-,p}
D_p B={B^{-1}}D_pB\end{equation*}
so we can retrieve the desired frequencies $\lambda_{k,p}=e^{\zeta_{k,p}}$ by computing $A_p$ and diagonalize it. Since all matrices $A_p$ are diagonalized by the same matrix $B$, the issue with grouping of the complex frequencies is easily solved just as in the previous case. The algorithm, which we call general domain ESPRIT, is summarized in Algorithm \ref{alg_gendomain}.
\begin{algorithm}
	\caption{General domain ESPRIT}\label{alg_gendomain}
	\begin{algorithmic}[1]
		\State Order the elements in $\Xi$ and $\Upsilon$.
        \State Form the general domain Hankel matrix $\bb{H}(\bb m,\bb n)=f( m+ n)$ from samples of $f$, where e.g. $\bb m$ is the order of the multi-index $m$ in $\Xi$.
		\State Compute the singular value decomposition $\bb{H} =\bb{U} \Sigma \bb{V}^*$.
		   \For {$p=1, \dots, d$}
		   \State Form $\bb{U}_{p+}$ and $\bb{U}_{p-}$ by deleting first and last elements in each fiber in the $p$:th coordinate, respectively.
		   \State Form $A_p =  ( \bb{U}_{p-}^\ast  \bb{U}_{p-})^{-1} \bb{U}_{p-}^\ast   \bb{U}_{p+}$.
		   \EndFor		
		\State Diagonalize $A_1=B^{-1} D_1 B$ by making an eigenvalue decomposition of $A_1$.
	    \For {$p=1, \dots, d$}
		   \State Compute the diagonal matrices $D_p = \dia(\lambda_{1,p}, \dots, \lambda_{K,p})=B A_p B^{-1}$.
		   \State Recover $\zeta_{k,p}=\log(\lambda_{k,p})$ (which are automatically correctly paired).
		\EndFor		
	\end{algorithmic}
\end{algorithm}

\section{Numerical examples}\label{numex}
To illustrate the methods discussed we will perform a number of numerical simulations. We have already shown Figure \ref{fig2dblock} and \ref{fig2dgen2}, in which 300 and 100 frequencies in 2d were retrieved up to machine precision, using Algorithm \ref{alg_dd} (in the first case) and Algorithm \ref{alg_gendomain} (in the second case). We continue here with a 3d example in the absence of noise, and then end by briefly discussing what happens in the presence of noise. Note that due to the non-linearity of the problem, this can behave quite differently in different situations. 

\begin{figure}
	\centering
	\subfloat[][Real part of a 3d data set sampled $21\times 21 \times 21$ lattice ($N=11$).]{\includegraphics[width=0.47\linewidth]{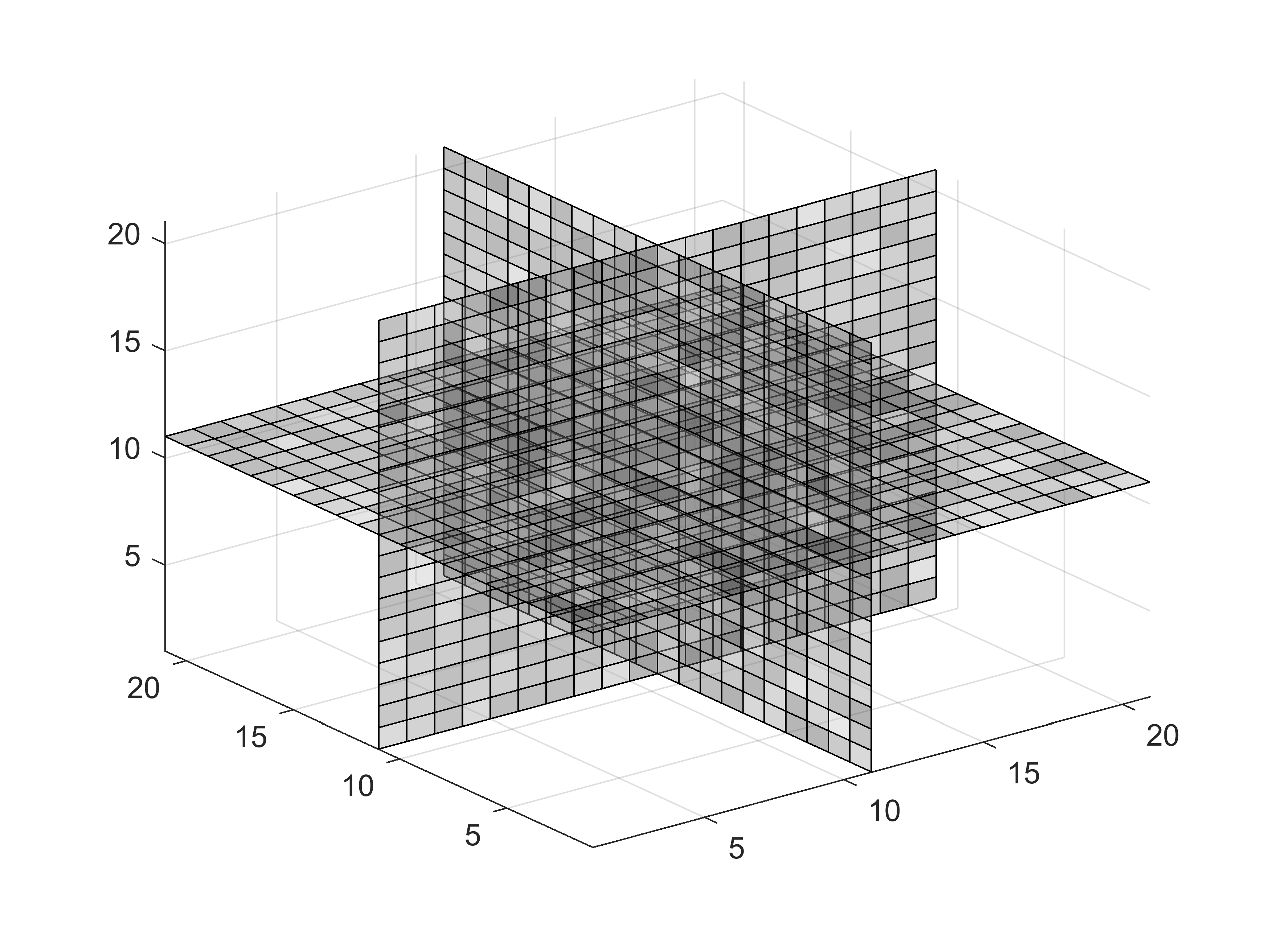}}
	\subfloat[][The corresponding block-Hankel matrix]{\includegraphics[width=0.47\linewidth]{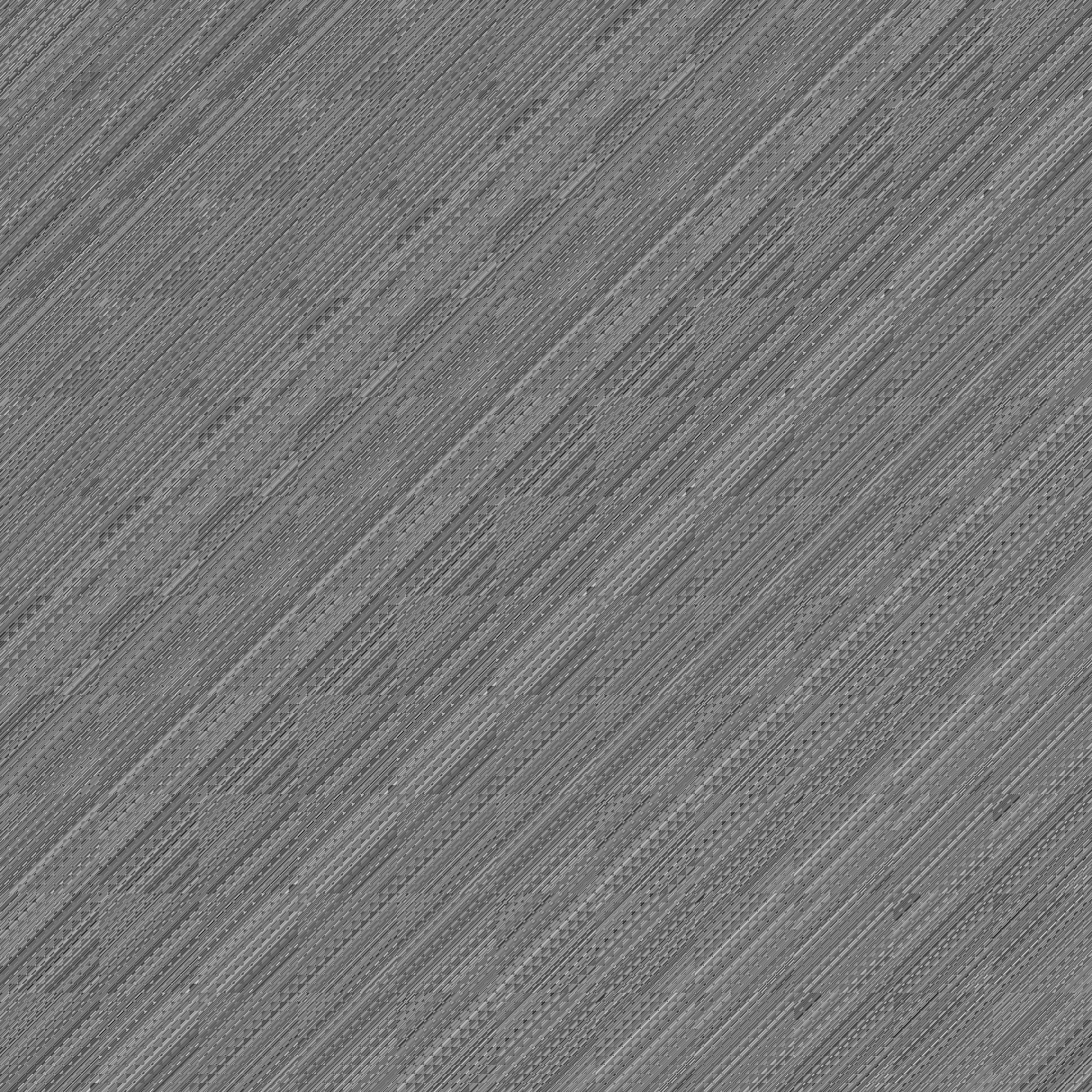}}\\
	\subfloat[][True 3d frequencies in blue dots and estimated frequencies in red circles]{\includegraphics[width=0.47\linewidth]{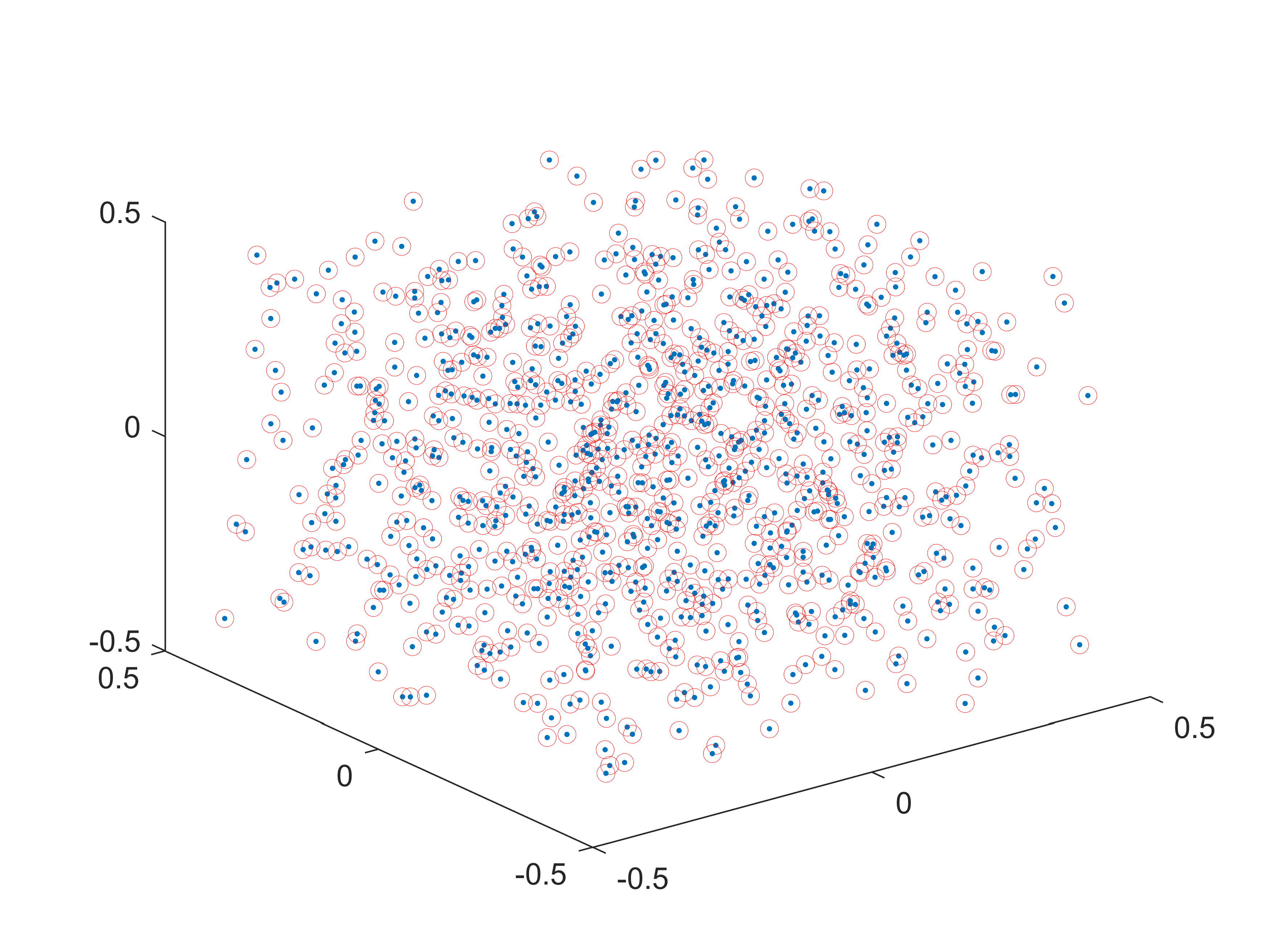}}
	\subfloat[][Error between estimated and true frequencies]{\includegraphics[width=0.47\linewidth]{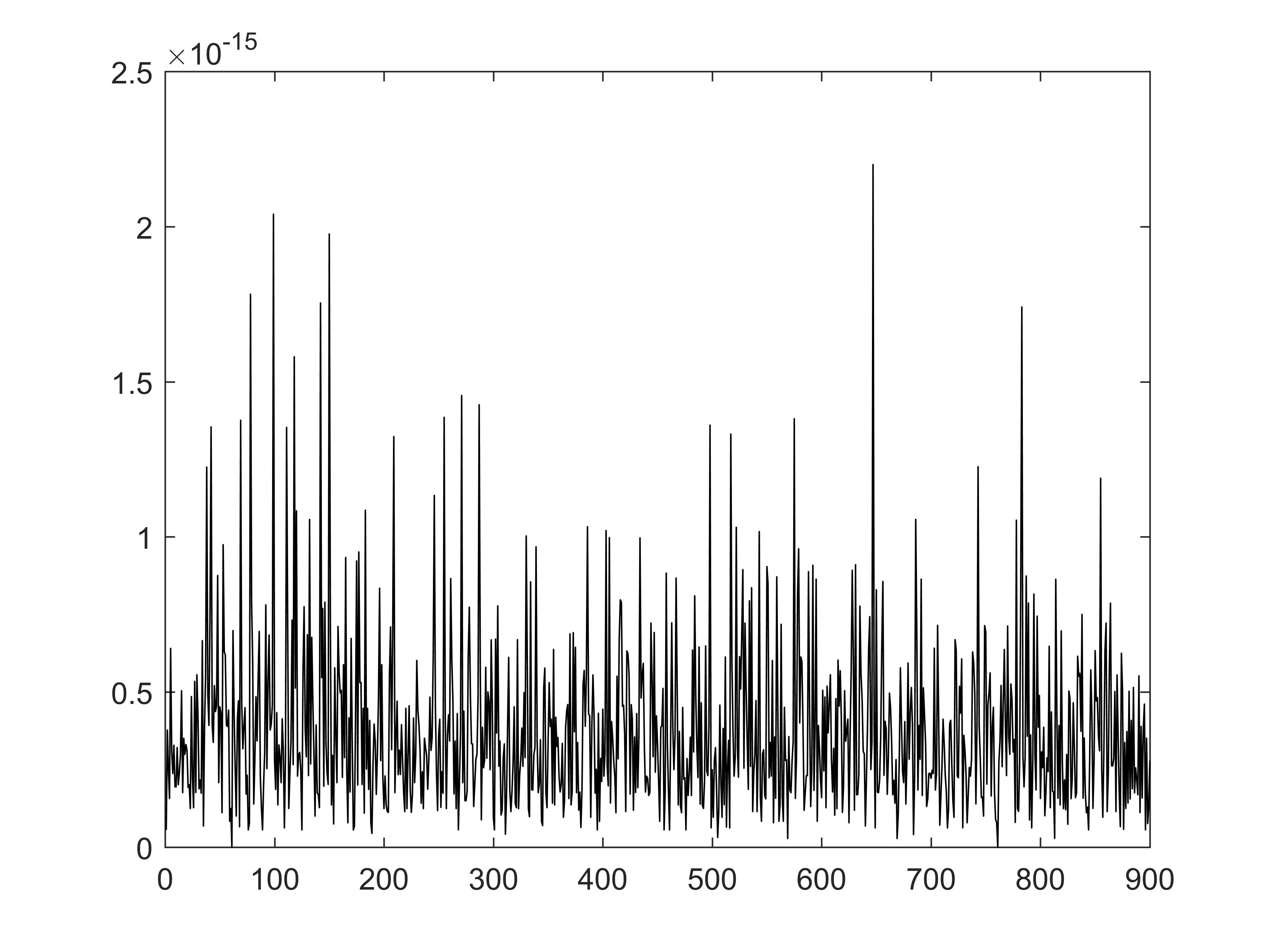}}
	\caption{\label{fig3dblock}Three dimensional example with data being a linear combination of 900 exponential functions.}
\end{figure}

The layout of Figure \ref{fig3dblock} is similar to that of Figure \ref{fig2dblock}. In this case 900 purely oscillatory exponential functions were used to generate a function $f$ sampled on a grid of size $21\times 21 \times 21$. In contrast to Figure \ref{fig2dblock}, the frequencies are now randomly distributed. The real part of $f$ is shown in panel a), the (real part of the) corresponding block-Hankel matrix is shown in panel b), the distribution of the exponentials in panel c) and the frequency reconstruction error is shown in panel d). According to \eqref{lok}, the maximal amount of frequencies we may retrieve in this situation is $11^2*10=1210$, to be compared with the total available data points $(21)^3=9261$.

Finally, we briefly illustrate the impact of noise on the proposed algorithms. However, we underline that the purpose of the present article is to provide an algorithm that correctly retrieves complex frequencies for functions of the form \eqref{lokf} in the absence of noise. If noise is present it is no longer true that the different $A_p$'s share eigenvectors, and hence choosing the eigenvectors of $A_1$ (step 8 in Algorithm \ref{alg_gendomain}) to diagonalize the rest becomes ad hoc. One way around this is to use some method to find the best simultaneous diagonalization of $\{A_p\}_{p=1}^d$ (see e.g. \cite{cardoso1996jacobi}), another is to use tools from optimization to preprocess the function $f$ so that it does become of the form \eqref{lokf} (see e.g.~\cite{andersson2017fixed}).

\begin{figure}
	\centering
	\includegraphics[width=0.9\linewidth]{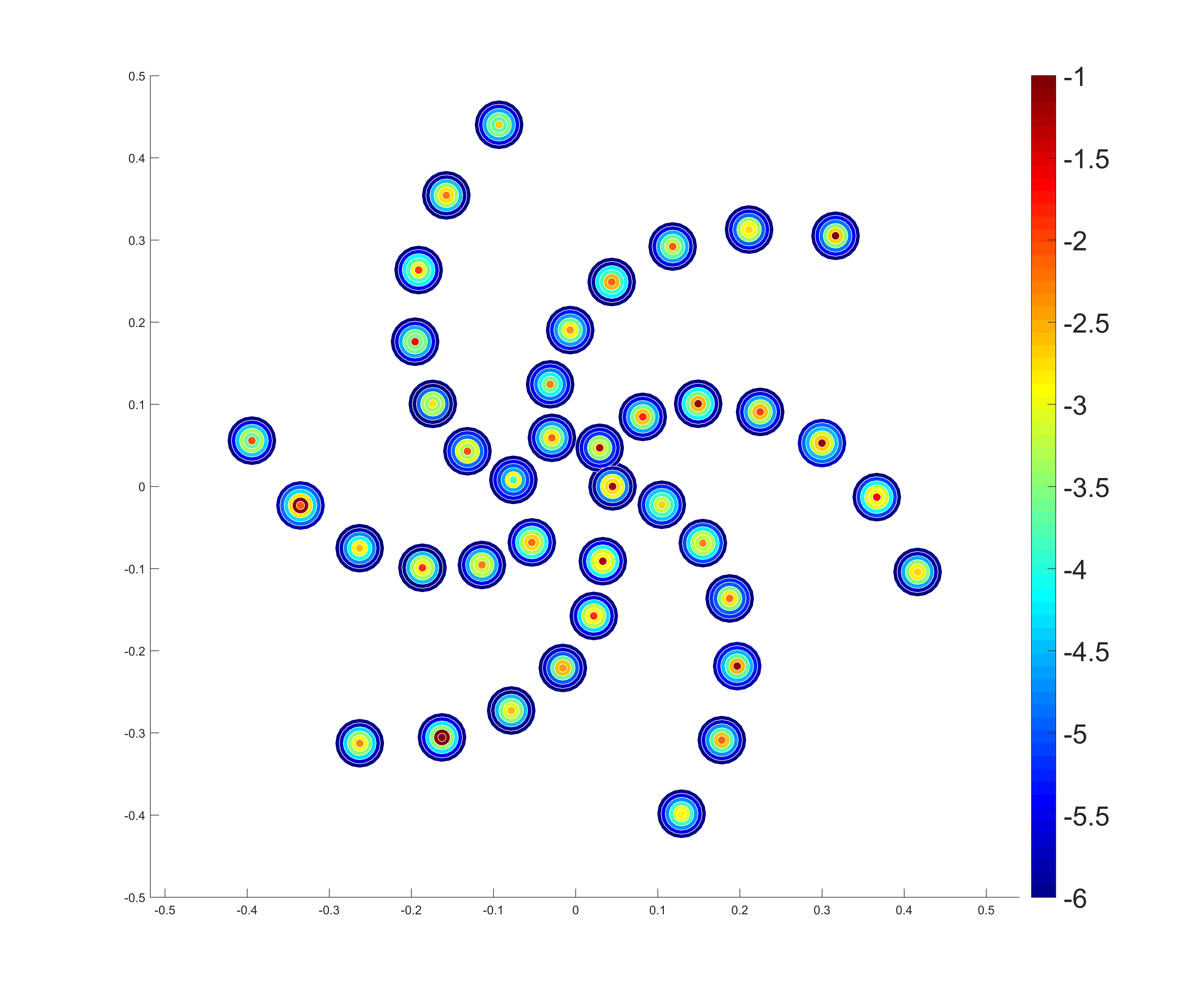}\label{<noise_nodes>}
	\caption{\label{fig2dnoise_a} Illustration of error in the estimation of frequencies. The size of the discs illustrate the estimation error at each frequency for the corresponding noise level in $\log_{10}$-scale.}
\end{figure}
\begin{figure}
	\centering
	\includegraphics[width=0.9\linewidth]{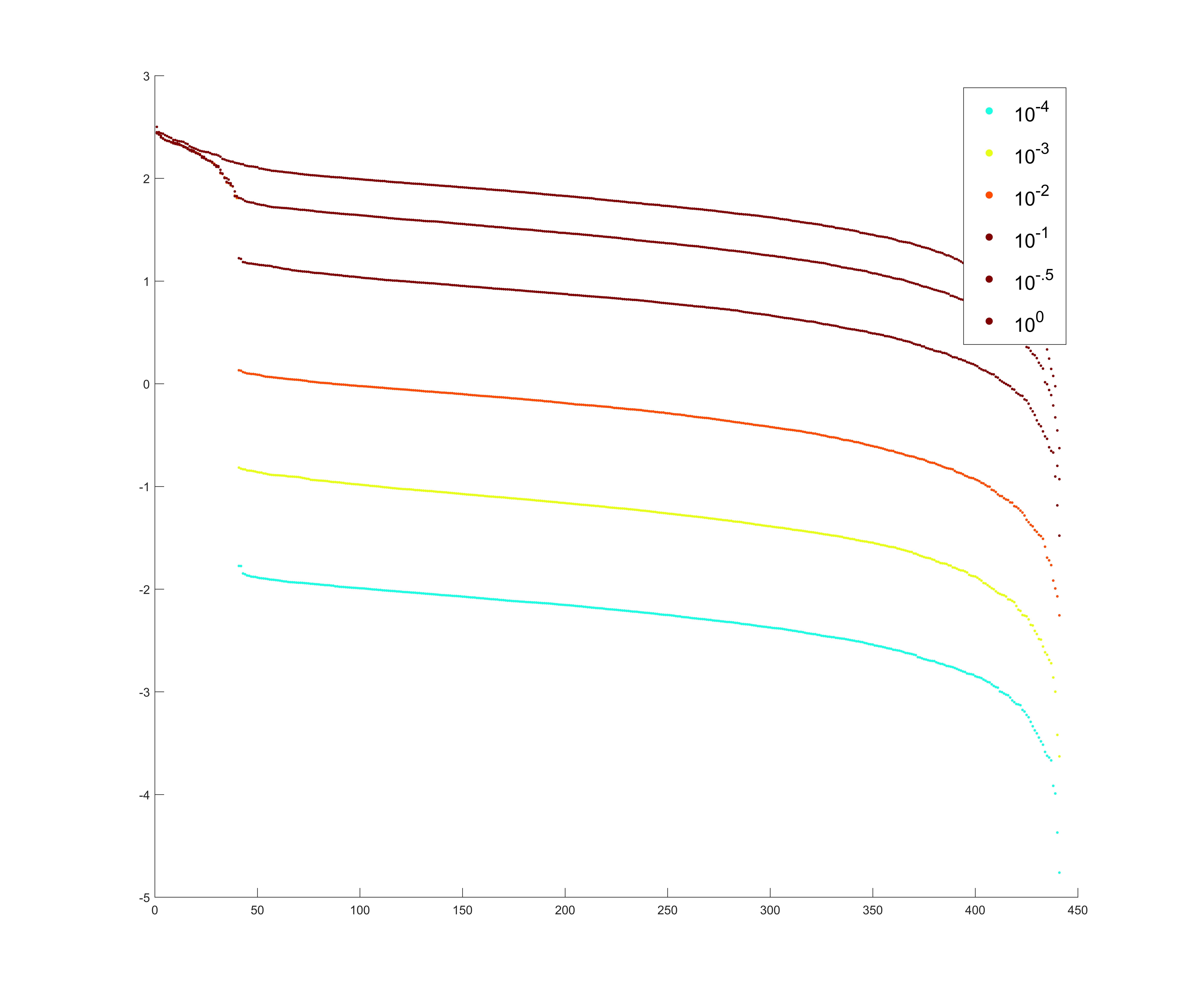}\label{<noise_sigma>}
	\caption{\label{fig2dnoise_b} The impact on the singular values for different noise levels.}
\end{figure}
The scatter plot in Figure \ref{fig2dnoise_a} shows the distribution of 40 frequencies in a case of purely oscillatory exponential functions. These are used to generate a function $f$ that is then sampled on a $41 \times 41$ grid. In addition, different levels of normally distributed noise is added to the original data. Six noise levels are chosen so that the ratio between the $\ell^2$ norm of the noise and the noise-free data is $10^0$, $10^{-.5}$, $10^{-1}$, $10^{-2}$, $10^{-3}$, and $10^{-4}$, respectively. Each noise level is portrayed with a particular color, according to the list in \ref{fig2dnoise_b}. The impact that the noise have on the singular values of the corresponding block-Hankel matrix is also shown in Figure \ref{fig2dnoise_b}. Recall that the noise free block Hankel matrix has rank 40, i.e. 40 non-zero singular values. When low noise is present this shows up as a jump in the magnitude of the singular values, clearly visible in plot \ref{fig2dnoise_b}. The first 40 singular values for e.g. yellow and orange are covered by the brown ones, therefore not visible. We can see that for the four lower levels, the impact of noise does not affect the original distribution on the singular values much, whereas for the highest noise level it affects essentially all singular values, and for the second highest level, it is just starting to have an impact. From this, we would expect to see a high error for the highest noise level, and a relatively small error for the four smallest noise levels.
To illustrate the effect on the individual frequency nodes, each of the nodes in Figure \ref{fig2dnoise_a}  have 6 circles around it. The color of each one of these circles show the error level in logarithmic scale according to the colorbar. Here we can see that the impact is pretty much as could be expected, with a low error for the lower noise levels and increasing for higher noise levels.

\section{Conclusions}
We have shown how to extract the underlying multi-dimensional frequencies from data that is constructed as linear combinations of (oscillating) exponential functions both for rectangular domains and more general domains. The approach does not require a pairing of one-dimensional frequencies components.

\bibliographystyle{plain}
\bibliography{MCref}

\end{document}